\documentclass[preprint,12pt]{elsarticle}

\usepackage{amsmath, amsthm,amssymb,bm, tikz, breqn, soul,setspace}
\usepackage{parskip, color}
\usepackage[left=1.2in,right=1.2in,top=1.2in,bottom=1.2in]{geometry}
\usepackage{mathtools}
\definecolor{mygray}{rgb}{0.95,0.95,0.95}
\usepackage{graphicx}
\usepackage{longtable}
\usepackage{hyperref}
\newcommand{\f}{{\mathbf f}}
\renewcommand{\a}{{\mathbf a}}

\newcommand{\R}{{\mathbb R}}
\newcommand{\E}{\mathcal{E}}

\newcommand{\hidederivation}[1]{}

\newtheorem{prop}{Proposition}
\theoremstyle{remark}
\newtheorem*{remark}{Remark}

\newcommand{\fnl}[1]{\big(#1\big)}

\journal{Computer Methods in Applied Mechanics and Engineering}

\begin{document}
\begin{frontmatter}

    \title{Computational Modeling of Coupled Interactions of Fluid Membranes with Embedded Filaments}

  \author[1]{Basant Lal Sharma}
  \address[1]{Department of Mechanical Engineering, Indian Institute of Technology Kanpur, Kanpur, 208016 UP, India}
 \author[2]{Luigi E. Perotti}
 \address[2]{Department of Mechanical and Aerospace Engineering, University of Central Florida, 12760 Pegasus Drive, Orlando, FL 32816}
   \author[3]{Sanjay Dharmavaram\corref{corr1}}
  \ead{sd045@bucknell.edu}
  \cortext[corr1]{Corresponding author}
 \address[3]{Mathematics Department, Bucknell University, Lewisburg, PA 17837}
  \begin{abstract}
    In this work, we present a computational formulation based on continuum mechanics to study the interaction of fluid membranes embedded with semiflexible filaments. This is motivated by systems in membrane biology, such as cytoskeletal networks and protein filaments aiding the cell fission process. We model the membrane as a fluid shell via the Helfrich-Canham energy and the filament as a one-dimensional Cosserat continuum. We assume the filament to be tethered to the surface of the membrane in a way that it is allowed to float on the surface freely. The novel filament-membrane coupling, which is anticipated to yield interesting physics, also gives rise to unique computational challenges, which we address in this work. We present validation results and apply the formulation to certain problems inspired by cellular biology. 
  \end{abstract}
\end{frontmatter}

\section*{Introduction}
\label{sec:introduction}

\begin{figure}[htb!]
    \centering
    \begin{tabular}{cc}
    \includegraphics[width=2.5in]{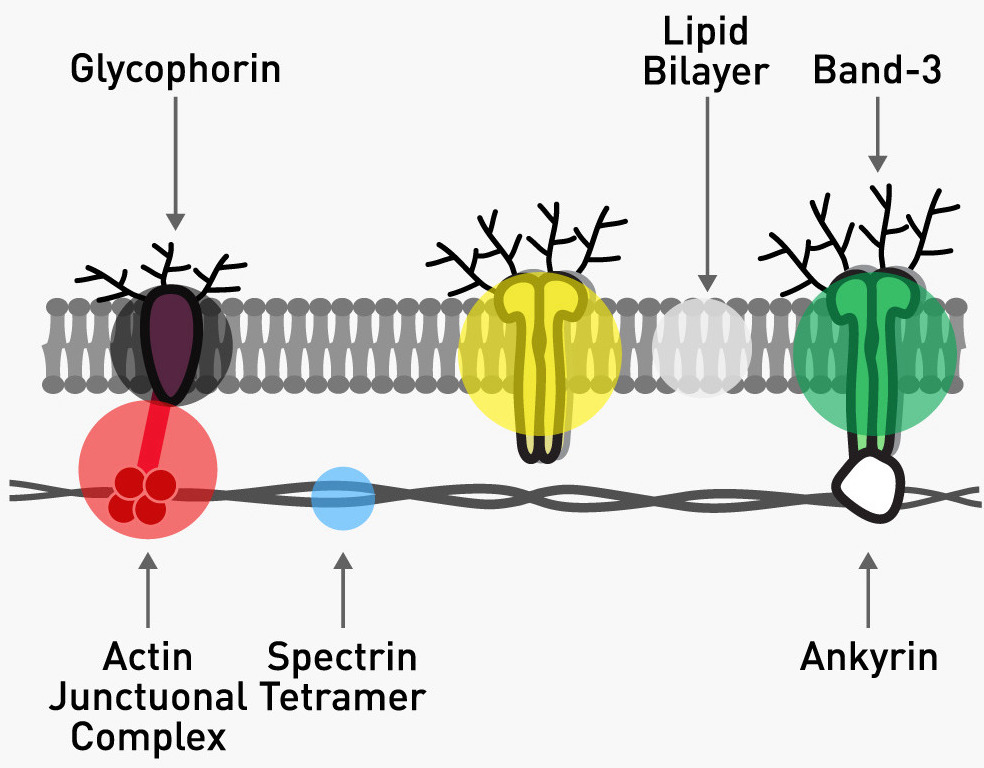} & 
    \includegraphics[width=3in]{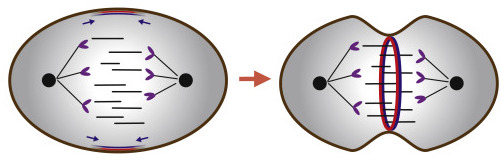}\\
    (a) & (b)
    \end{tabular}
    \caption{(a) Schematic showing actin filament anchored on the membrane through membrane-bound proteins, courtesy Fig.~1B from \cite{tang2017openrbc}. (b) Schematic showing a contractile actomyosin forming along the mid-plane, courtesy Fig.~3B \cite{cheffings2016actomyosin}.}
    \label{fig:actin}
\end{figure}

Many cellular processes involve the interaction between cell membranes with filamentous structures. For instance, most cells contain a \emph{cytoskeleton} which is a dynamic network of protein filaments crisscrossing the cell's cytoplasm and the cell membrane. The cytoskeleton not only provides structural rigidity to the cell but also plays a vital role in the cell's signaling pathways through the transport of cargo-filled vesicles.  In most eukaryotic cells, a vital component of the cytoskeleton is the \emph{cortex}, a network of semiflexible actin filaments coating the underside of the cell membrane \cite{chalut2016actin}, which gives cells their characteristic shape and rigidity. It also plays a vital role in cell mobility and cell mechanics. Fig.~\ref{fig:actin}(a) (courtesy \cite{tang2017openrbc}) shows a schematic of actin filaments tethered to the membrane with membrane-bound proteins, such as \emph{glycophorin} and \emph{band-3}-ankyrin. 
The motivation for the present work comes from this as first example of such structures found in cellular biology.  

A second  example, where the interaction of the cell membrane with embedded filaments is relevant, refers to the events during the formation of lipid vesicles \cite{guizetti2011cortical, akamatsu2020principles}. The latter play an important role in the transport and sequestration of biochemical cargo from cells. Vesicle formation is typically initiated by a collection of protein complexes, collectively called the ESCRT-III (Endosomal Sorting Complex Required for Transport) machinery. Key components of this complex are protein monomers that adhere to the cell membrane and polymerize into filaments (approximately 17nm thick) \cite{guizetti2011cortical}. The filaments proceed to form a helical constricting bundle, which causes the formation of a bud and eventually cleaves it into a spherical vesicle. The exact mechanism by which the helical filament induces curvature in the membrane is poorly understood \cite{mccullough2020membrane}. Some studies have shown that the filaments form a flat spiral on lipid membranes and suddenly undergo a buckling instability into a helical state \cite{mccullough2020membrane,chiaruttini2015relaxation}. The present work is the first step in this direction that prescribes one way to model this.

A third example to motivate our work is the mechanics of cell division aided by the formation of an \emph{actomyosin ring} \cite{cheffings2016actomyosin}. During \emph{cytokinesis} (the cell fission process), a contractile ring made up of polymerized protein complexes (composed of actin and myosin II) forms at the equatorial position between the two daughter cells, shown schematically in Fig.~\ref{fig:actin}(b) (courtesy \cite{cheffings2016actomyosin}). Aided by ATP, the myosin motor proteins move along the actin filaments and begin constricting the ring \cite{o2018mechanisms}. Understanding the mechanics of cytokinesis has implications in synthetic biology -- an emerging interdisciplinary area encompassing biology, material science, and engineering aimed to design synthetic biological systems. The challenges encountered in replicating cytokinesis in a synthetic cell are discussed in a recent review article \cite{baldauf2022actomyosin}. In particular, it is noted that in such a system, anchoring the actomyosin ring along the equatorial position of the cell is not an easy task, as upon myosin activation, the ring slips along the membrane. Understanding the mechanics of membrane-filament interaction can potentially shed light on the mechanisms for controlling this process. 

A common theme in the three examples, briefly described above, is the mechanical coupling between the lipid membrane and the embedded protein filaments.  The goal of this work is to propose a computational framework based on continuum mechanics for understanding the coupled interaction of fluid membranes with embedded semiflexible filaments.

\subsection*{Structure of Lipid Membranes and Protein Filaments}

Lipid membranes are made up of lipid molecules arranged in a bilayer structure. The hydrophobic tails of the molecules are embedded on the inner side of the bilayer, and the polar heads are exposed to the aqueous cytoplasm. Within each layer, the molecules can freely move, thus rendering a fluid-like behavior to the membrane, similar to a soap film. Unlike the latter, the membrane resists bending as doing so exposes the hydrophobic tails to the aqueous environment of the cytoplasm. This ``in-plane'' fluidity of the membrane ascribes unique mechanical properties to the membrane. Often, membrane-bound proteins float on the surface and play a key role in various cellular processes.

Most biological filaments at the cellular level are polymers of proteins. They can range in diameter from 7nm (e.g., actin filaments) to 25nm (e.g., microtubules). Protein filaments are critical to the proper functioning of the cell. They perform a variety of roles, from providing structural support to acting as highways for material transport. Since protein monomers that constitute a filament are chiral molecules, the filament inherits an overall twisted chiral structure \cite{holmes1990atomic}. From a mechanistic viewpoint, they behave as semiflexible filaments which resist extension, bending, and twisting. The Cosserat rod model is a versatile framework for a continuum description of filaments \cite{floyd2022stretching}. Twist-bend coupling in actin filaments \cite{enrique2010origin,bibeau2023twist} and chirality of the filament can be easily incorporated within this description; however, in the present work, we do not investigate the latter aspect of filaments for simplicity.

As a consequence of the in-plane fluidity of the membrane, the filaments may be considered to be floating on the surface. The schematic shown in Fig.~\ref{fig:actin}(a) illustrates this idea. Even though the filaments are anchored on the membrane with the membrane-bound proteins, the latter are freely floating on the surface. In fact, it is a unique feature of the model systems considered in this work that the filaments can freely drift on the surface without encountering any resistance.

\subsection*{State-of-the-art in Modeling}
Current approaches in modeling this class of systems can be broadly classified into two types: (1)  molecular dynamics simulations and (2) continuum models. 

Molecular dynamics studies such as \cite{mandal2020molecular} have used fully atomistic simulations to model the ESCRT-membrane interactions, but such models are computationally extremely intensive, and only relatively small systems can be simulated. A slightly more coarse-grained approach where the filament is modeled as a set of beads connected by strings has been proposed in \cite{harker2019changes}.
A two-component coarse-grained molecular dynamic approach to study the coupled mechanics of membrane and its cytoskeletal network in red blood cells has been proposed in \cite{li2012two}. In this method, the lipid molecules, actin proteins and band-3 anchor proteins (c.f., Fig.~\ref{fig:actin}(a)) are each modeled as a ``particle'' with specific position and orientation-dependent interactions. The drawback of this method is the preponderance of phenomenological constants, as fitting parameters, for the interactions in order to get meaningful results.

Continuum models so far have been limited to modeling cytoskeletal networks. In the works of Boey, Discher, et al. \cite{boey1998simulations,discher1998simulations}, the cytoskeletal network of erythrocytes (red blood cells) is modeled as a network of Hookean springs. While their work only considered flat membranes and neglected the underlying lipid membrane, a more realistic mechanical description of erythrocytes has been considered in \cite{peng2013lipid, rangamani2014multiscale}. In this latter study, the lipid membrane and the cytoskeletal network are each represented by a triangular mesh. Equilibrium configurations are determined using dissipative particle dynamics. A crucial modeling assumption in this work is that the nodes of two meshes are identified. While this may be sufficient for cytoskeletal networks, it is not straightforward to generalize the methodology to the other biological systems noted above.

In this work, we explore the mechanics of membrane and filament coupling through a continuum model for both structures. 

\subsection*{Outline}
This paper is organized as follows. In Sec.~\ref{sec:formulation} we present the mathematical formulation for the system--a fluid membrane embedded with a filament. We follow a variational approach and postulate an energy functional that combines the Helfrich-Canhman model \cite{helfrich1973elastic, jenkinsRBC} for the lipid membrane with a Cosserat rod model \cite{Cosserat,antman} for the protein filament. The fact that the filament is embedded, yet floating on the fluid membrane, causes computational challenges which we discuss and address in Sec.~\ref{sec: computational challenges}. We present the weak form of the energy functional in Sec.~\ref{sec:weak form}. In Sec.~\ref{sec:numerical results} we compute equilibrium configurations of the system by discretizing the weak form of the functional. We use the Loop subdivision finite element methods \cite{loop1987smooth, cirak2000subdivision} to discretize the membrane and Hermite cubic finite elements for the filament. We present validation studies to test the accuracy of our model and consider some applications. We conclude by pointing out the prospects and limitations of our approach in Sec.~\ref{sec:discussion and conclusions}.

\section*{Notation}
{We follow the standard notation in mathematics and mechanics \cite{gurtin1982introduction}.} The set of real numbers is denoted by $\mathbb{R}$ whereas the three dimensional Euclidean space is denoted by $\mathbb{R}^3$, which is equipped with the standard inner-product {denoted with a dot, i.e., $\mathbf{a}\cdot\mathbf{b}$ for given $\mathbf{a},\mathbf{b}\in\mathbb{R}^3$}. 
We use boldface font to represent the vectors and tensors. 
The second order tensor $\mathbf{a}\otimes\mathbf{b}$, for given $\mathbf{a},\mathbf{b}\in\mathbb{R}^3$, defined as a linear transformation from $\mathbb{R}^3$ to $\mathbb{R}^3$ such that $(\mathbf{a}\otimes\mathbf{b})[\mathbf{x}]=\mathbf{a}(\mathbf{b}\cdot\mathbf{x})$ for all $\mathbf{x}\in\mathbb{R}^3$.
We denote the transpose of a tensor in a way that $\mathbf{A}^T$ represents the transpose of $\mathbf{A}$ (defined so that $(\mathbf{a}\otimes\mathbf{b})^T=\mathbf{b}\otimes\mathbf{a}$).
Sometimes for brevity, and for emphasis on other occasions, we abuse the notation and alternate between a comma notation and an explicit form, i.e., $
f_{,x}(x,y)=\partial_x f(x,y)=\frac{\partial}{\partial x} f(x,y)$, $f_{,xy}= \partial_{xy}f(x,y)=\frac{\partial^2}{\partial y\partial x} f(x,y)$, etc.
The square brackets, in general, denote the linear action of an operator (say the gradient of a function, the first variation of functional, etc), which is mostly clear from the context with an exception of the use of similar notation for a closed interval of the real line.

\section{Formulation of the Model}
\label{sec:formulation}
\begin{figure}[htb!]
    \centering
    \begin{tabular}{cc}
    \includegraphics[width=2in]{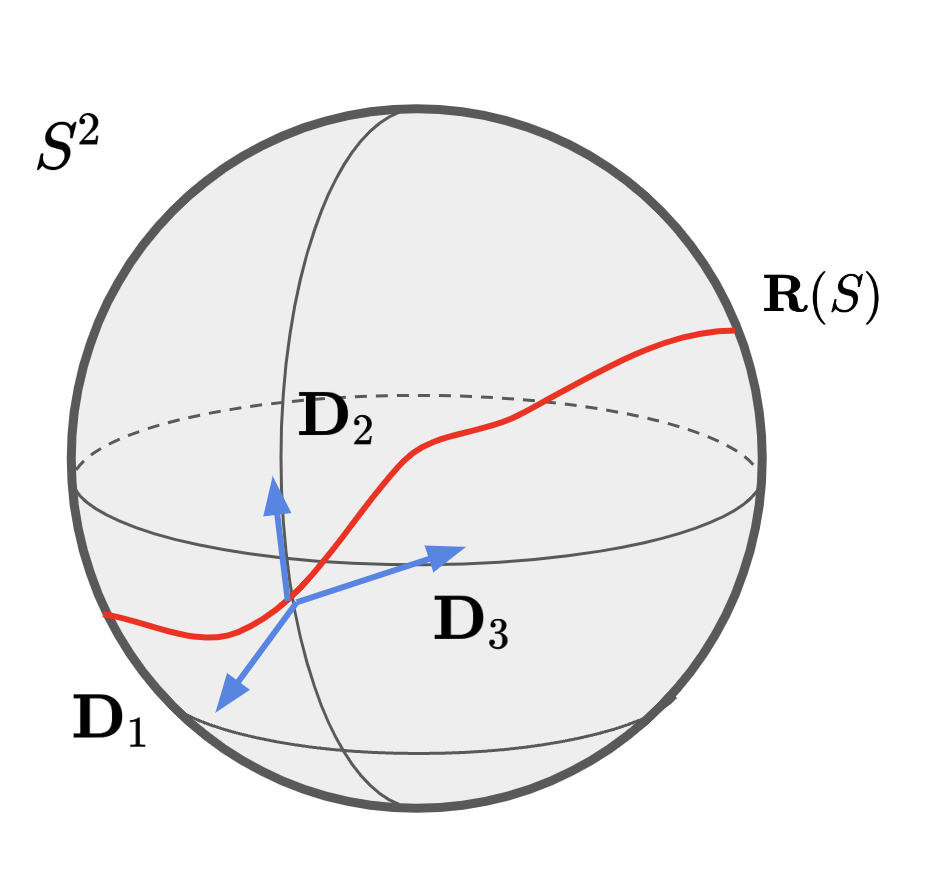} & 
    \includegraphics[width=2in]{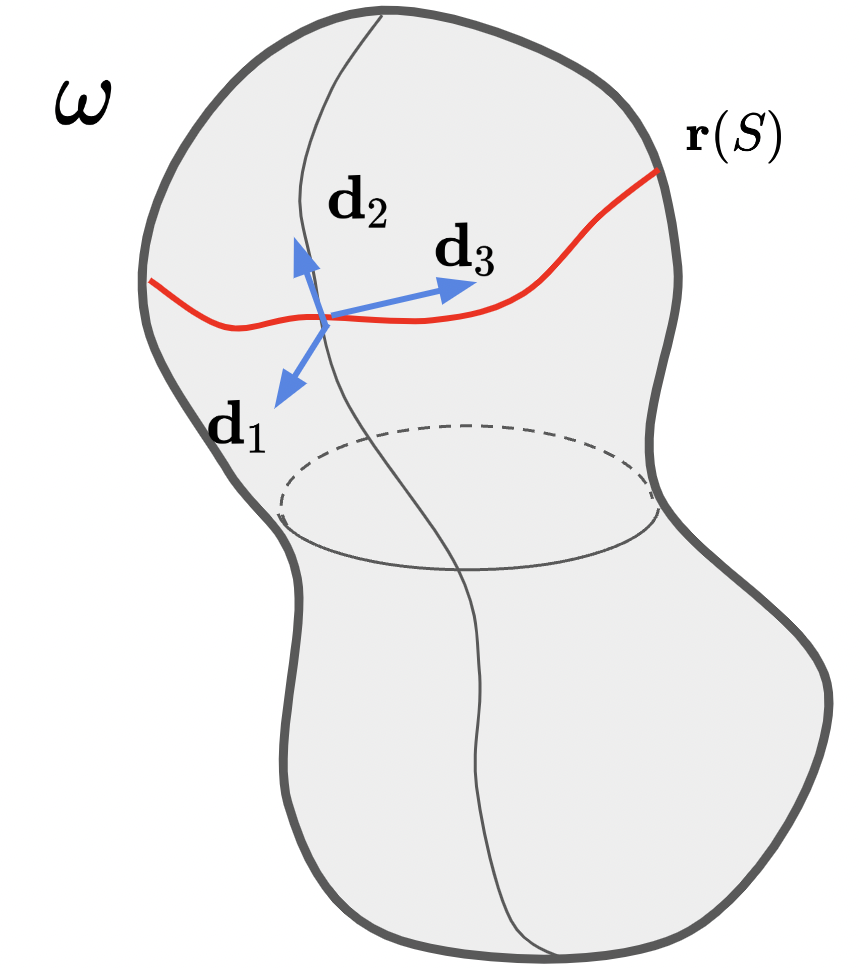}\\
    (a) & (b)
    \end{tabular}
    \caption{Schematic showing (a) reference and (b) deformed configurations of the membrane-rod system. Director in the reference and deformed configuration are shown as arrows.}
\label{fig:schematic_formulation}
\end{figure}

In this work, we model the membrane as a {\em fluid shell} with a spherical reference configuration. While working with appropriate dimensional variables, we assume that the reference configuration of the membrane is a unit sphere $S^2$. The filament is modeled as a \emph{special Cosserat rod} which is assumed to be embedded on the surface of the membrane. Due to the in-plane fluidity of the membrane, it is a natural assumption that the filament is allowed to float freely on the surface; in particular, the filament applies reaction forces only in the direction normal to the surface.  Our focus in this work is the computation of various equilibrium configurations (statics) of such a membrane-filament system, and, therefore, we do not consider any viscous forces exerted by the membrane on the filament. A schematic for the fluid shell-filament system is shown in Fig.~\ref{fig:schematic_formulation}.  

\subsection{Lipid membrane: Helfrich-Canham Model}
\label{ssec:Helfrich}
We model the membranes' stored energy due to elastic bending via the Helfrich-Canham energy functional \cite{helfrich1973elastic, jenkinsRBC}. An important feature of this form of energy is that it naturally incorporates the in-plane fluidity of membranes. From a microscopic viewpoint, the fluidity in lipid membranes arises because of the disordered arrangement of the hydrophobic tails of the lipid molecules. The lipids molecules are free to move around along the surface of the membrane and experience very little shear force in this process. 

A consequence of in-plane fluidity of the membrane is the lack of a well-defined reference configuration for the membrane. In this work, we use the term \emph{reference configuration} of the shell to denote any convenient parametrization of the surface. We denote 
\[
\f:S^2\rightarrow \omega \subset \R^3
\]
to be the deformation map for \emph{current configuration} of the shell (the deformed surface $\omega$). Schematics for the reference and deformed configuration are shown in Fig.~\ref{fig:schematic_formulation}(a) and (b). Let $\mathbf{X}:=(X^1,X^2)\in\R^2$ denote a chosen (local) parametrization for the reference unit sphere $S^2$. The Helfrich-Canham energy functional \cite{canham1970minimum, helfrich1973elastic} is given by:
\begin{equation}
  \E_{\text{Helf}}\fnl{\mathbf{f}}=\int_\omega\left( c (H-C_0)^2 + c_g K \right)\;da,
  \label{eq:energy_pre}
\end{equation}
where $H$ and $K$ are the mean and Gaussian curvatures of $\omega$, while $c$ and $c_g$ are their respective bending stiffnesses; also, $C_0$ is the preferred curvature of the membrane, and $da$ is the area measure on $\omega$. The explicit dependence of the energy on the deformation map $\mathbf{f}$ is shown by the term in parentheses on the left side of the above equation. For closed surfaces (i.e., without a boundary), by virtue of the Gauss-Bonnet theorem \cite{carmo1992riemannian}, the Gaussian curvature term (the second term in \eqref{eq:energy_pre}) integrates to a topological constant. Thus, this term does not play a role in our work as we only consider topologically same (spherical) configurations. Additionally, since the area modulus of lipid membranes is significantly larger than the bending moduli,  we assume area incompressibility and impose the following total area constraint:
\begin{equation}
  \int_\omega da = 4\pi,
  \label{eq:constr}
\end{equation}
where the right side of the equation is the surface area of the reference spherical state (the unit sphere $S^2$). Note that  \eqref{eq:constr} is a \emph{global} constraint on the surface area of $\omega$. While it is conventional in continuum mechanics to model incompressibility by fixing the area strain at every point on the surface to unity \cite{jenkinsRBC}, such a \emph{local} incompressibility constraint is not convenient for numerical computations, especially with the finite element methods, as it leads to a mixed formulation. For topologically spherical membranes (the only ones considered in this work), the two formulations have been shown to be equivalent \cite{dharmavaram2015equivalence} as far as equilibria and their stability is concerned.

We assume further that the membrane supports an internal osmotic pressure $p$, so that the total free energy of the membrane is given by 
\begin{equation}
  \E_{m}\fnl{\mathbf{f}} = \mathcal{E}_{\text{Helf}}\fnl{\mathbf{f}}-pV,
  \label{eq:free energy membrane}
\end{equation}
where $V$ is the volume enclosed by $\omega$ (cf., \eqref{eq:energy_pre}). 

\subsection{Filament: Special Cosserat Rod}
\label{ssec:cosserat}
We model the rod (a term which we use interchangeably with `filament') using the framework of special theory of Cosserat  \cite{antman} in which the rod's center curve is described by $\mathbf{r}(S)$. Here, $S\in[0,L]$ denotes the (referential) arclength parameter with $L$  being the reference length of the rod (measured in dimensions of the radius of the membrane). Thus, 
\[
\mathbf{r}: [0,L]\to\mathbb{R}^3.
\]
The rod cross-section is spanned by two directors $\mathbf{d}_1(S)$ and $\mathbf{d}_2(S)$, and we assume that there are three (orthonormal and oriented) directors $\{\mathbf{d}_1(S),\mathbf{d}_2(S),\mathbf{d}_3(S)\}$. In the reference configuration, the center curve and the directors are described by $\mathbf{R}(S)$, $\{\mathbf{D}_1(S),\mathbf{D}_2(S),\mathbf{D}_3(S)\}$ (orthonormal and oriented). Note that $S$ satisfies the condition 
\[
\lVert\mathbf{R}'(S)\rVert\equiv 1.
\]
In Fig.~\ref{fig:schematic_formulation}, the directors in the reference and deformed configurations are represented by arrows. The directors, $\mathbf{d}_i,\;i\in\{1,2,3\}$, form an orthonormal and oriented basis (same statement for $\mathbf{D}_i,\;i\in\{1,2,3\}$); in particular, $\mathbf{d}_1\cdot(\mathbf{d}_2\times\mathbf{d}_3)=+1$. This assumption is the characteristic of the \emph{special Cosserat rod model} and it is a reasonable assumption for  thin biological filaments where cross-sectional shear can be neglected. Thus, we can map the current directors to reference directors (both parametrized by $S$), in other words associate the cross-section in the deformed configuration,  via a rotation tensor 
\[
\mathbf{Q}(S)=\sum_{i=1}^3\mathbf{d}_i(S)\otimes\mathbf{D}_i(S).
\]
Note that $\mathbf{Q}:\mathbb [0,L]\to \text{Orth}^+$. Here $\text{Orth}^+$ denotes the set of rotations in three dimensions.
That is,
\begin{equation}
\mathbf{d}_i(S) = \mathbf{Q}(S)\mathbf{D}_i(S),\;\text{ for }i=1,2,3 \text{ and }S\in[0,L].
\label{eq:d=QD}
\end{equation}
The strain measures for the rod are given by the derivatives of $\mathbf{r}$ and $\mathbf{Q}$. We define
\begin{equation}
\nu_i(S) := \mathbf{r}'(S)\cdot\mathbf{d}_i(S),\;i=1,2,3 \text{ for }S\in[0,L],
\label{eq:nu_i}
\end{equation}
where $\nu_1(S)$ and $\nu_2(S)$ are the shear strains at the cross-section $S$, and $\nu_3(S)$ the stretch along the rod at $S$. 
The bending and twisting strains in the rod are related to the derivative of the directors:
\begin{equation}
\mathbf{d}_i'(S) = \bm{\kappa}(S)\times \mathbf{d}_i(S) + \mathbf{Q}(S)\mathbf{D}_i'(S),\;i=1,2,3\text{ for }S\in[0,L],
\label{eq:d'}
\end{equation}
where $\bm{\kappa}$ is an axial vector defined by $\bm{\kappa} := \text{axial}(\mathbf{Q}'\mathbf{Q}^T)=\sum_{i=1}^3\kappa_i\mathbf{d}_i$.
In general, prime refers to differentiation with respect to $S$.
The quantities $\kappa_1$ and $\kappa_2$ can be interpreted as bending strains in the rod (of the center curve), while $\kappa_3$ as the twisting strain of the cross-section about the center curve. 

Note that the second term on the right-hand side of \eqref{eq:d'} is relevant if the reference configuration of the rod is curved. Since $\mathbf{D}_1$, $\mathbf{D}_2$, $\mathbf{D}_3$ form an orthonormal basis of $\mathbb{R}^3$, we may assume without loss of generality that
\begin{equation}
    \mathbf{D}'_i(S) = \sum_{j=1}^3 \alpha_{ij}(S)\mathbf{D}_j(S),\;i=1,2,3,
    \label{eq:D'}
\end{equation}
for prescribed functions $\alpha_{ij}$, which depend on the curvatures of the reference configuration of the rod. Orthonormality of the directors, naturally, implies that $\alpha_{ij}=-\alpha_{ji}$, and in particular, $\alpha_{ii}=0$, for $i,j\in\{1,2,3\}$. 
In this work, we treat $\alpha_{ij}$ as constants; in this case, these parameters can be related to the Frenet–Serret parameters of the reference curve. For example, for a rod with a circular reference state, without loss of generality, we assume that its reference state lies on the $xy$ plane,  that is
\begin{equation}
    \mathbf{R}(S) = \frac{L}{2\pi} ( \mathbf{i}\cos S +\mathbf{j} \sin S),
\end{equation}
with
\begin{subequations}    
\begin{eqnarray}
    \mathbf{D}_1(S) &=& ( \mathbf{i}\cos S + \mathbf{j}\sin S),\\
    \mathbf{D}_2(S) &=& \mathbf{k},\\
    \mathbf{D}_3(S) &=& (-\mathbf{i} \sin S + \mathbf{j}\cos S).
\end{eqnarray}
\label{eq:D1 D2 D3}
\end{subequations}
Thus,
$    \mathbf{D}'_1(S) = \mathbf{D}_3,\;\mathbf{D}_2'(S)=0,\;\mathbf{D}_3'=-\mathbf{D}_1,$
which yields $\alpha_{13}=-\alpha_{31}=1$, $\alpha_{12}=\alpha_{21}=\alpha_{32}=\alpha_{23}=0$.

Plugging \eqref{eq:D'} into \eqref{eq:d'}, and using \eqref{eq:d=QD}, we obtain 
\begin{equation}
    \mathbf{d}_i' = \bm{\kappa}\times\mathbf{d}_i + \sum_{j=1}^3\alpha_{ij}\mathbf{d}_j,\;i=1,2,3.
\end{equation}
Using the orthonormality of the directors, the previous equation gives us the following equations for the bending and twisting strains:
\begin{subequations}
\begin{eqnarray}
    \kappa_1 = \mathbf{d}_2'\cdot\mathbf{d}_3-\alpha_{23} = -(\mathbf{d}_3'\cdot\mathbf{d}_2-\alpha_{32}),
    \label{eq:kappa1}\\
    \kappa_2 = \mathbf{d}_3'\cdot\mathbf{d}_1-\alpha_{31} = -(\mathbf{d}_1'\cdot\mathbf{d}_3-\alpha_{13}),
    \label{eq:kappa2}\\
    \kappa_3 = \mathbf{d}_1'\cdot\mathbf{d}_2-\alpha_{12} = -(\mathbf{d}_2'\cdot\mathbf{d}_1-\alpha_{21}).
    \label{eq:kappa3}
\end{eqnarray}
\label{eqs:kappa_i}
\end{subequations}
In the framework of hyperelastic materials, the strain variables \eqref{eq:nu_i} and \eqref{eqs:kappa_i} can be used to postulate a general constitutive law \cite{antman},
\begin{equation}
\mathcal{E}_{rod}\fnl{\mathbf{r},\mathbf{Q}}=\int_{0}^L W(\nu_1,\nu_2,\nu_3,\kappa_1,\kappa_2,\kappa_3)\;dS.
    \label{eq:E rod W}
\end{equation}
The explicit form for $W$ depends on considerations of material symmetry for the rod. Note that in \eqref{eq:E rod W}, the dependence of $\mathcal{E}_{rod}$ on $\mathbf{r}$ and $\mathbf{Q}$ reflects the dependence of energy on the two kinematic variables of the rod model.

In this work, we assume that the rod is unshearable which is a reasonable assumption for thin filaments. Accordingly, we set
\begin{equation}
    \nu_1=\nu_2\equiv 0.
    \label{eq:unshearable}
\end{equation}
To keep the discussion simple, we assume an isotropic and compressible constitutive law for the rod with the following form:
\begin{equation}
    \mathcal{E}_{rod}=\int_0^L \frac{D}{2}(\nu_3-1)^2 + \frac{E}{2}(\kappa_1^2+\kappa_2^2+\kappa_3^2)\;dS,
    \label{eq:E rod isotropic}
\end{equation}
where $D$ is the elastic modulus of the rod, and $E$ is the bending modulus. Note that $E$ can also be interpreted as the twisting modulus as it is also the coefficient of the twisting strain $\kappa_3$. Even though we restrict ourselves to an isotropic constitutive law, our formulation is general enough to accommodate more realistic models for chiral filaments (e.g., hemitropic constitutive law \cite{healey2017symmetry}). While it may seem necessary to enforce the constraints \eqref{eq:unshearable} either using Lagrange multipliers or a penalty term, we find that this is not needed in the scheme we describe below.

\section{Computational Challenges}
\label{sec: computational challenges}
Our goal is to compute the equilibrium configurations associated with the energy (cf., \eqref{eq:free energy membrane} and \eqref{eq:E rod isotropic}):
\begin{equation}
    \mathcal{E}\fnl{\mathbf{f},\mathbf{r},\mathbf{Q}} =\mathcal{E}_{m}\fnl{\mathbf{f}}+\mathcal{E}_{rod}\fnl{\mathbf{r},\mathbf{Q}},
    \label{eq:energy 1}
\end{equation}
subject to the constraints \eqref{eq:constr} and \eqref{eq:unshearable}.  For emphasis, in the above equation, we show the explicit dependence on the various degrees-of-freedom ($\mathbf{f}$, $\mathbf{r}$, $\mathbf{Q}$) in parentheses. Recall that due to \eqref{eq:d=QD}, $\mathbf{Q}$ can be interchanged with the three directors $\mathbf{d}_i$, $i=1,2,3$.

We encounter two main computational challenges as we try to find the critical points of \eqref{eq:energy 1}.
\begin{enumerate}
    \item \textbf{Membrane-rod coupling}: Since the rod remains tethered to the surface of the membrane, enforcing this embedding condition naively entails imposing a point-wise constraint:
\begin{equation}
    \lVert\mathbf{r}(S)-\mathbf{f}(X^1,X^2)\rVert=0,
    \label{eq:naive constraint}
\end{equation}
where the norm is the Euclidean norm in $\mathbb{R}^3$; thus there exits a tuplet $(X^1,X^2)$ corresponding to given $S$ such that this condition holds. Point-wise constraints such as the one above are computationally cumbersome to impose. 
  \item \textbf{Membrane fluidity:} Applying mesh-based computational methods to fluid membranes is often plagued with extreme mesh distortion and the presence of zero-energy modes \cite{ma2008viscous, feng2006finite}. These stem from the in-plane fluidity of the constitutive law. A variety of numerical schemes have been proposed to address this issue \cite{feng2006finite, ma2008viscous, Elliot, zhao2017direct, healey2017symmetry} which either add in-plane viscosity (there by changing the constitutive law) or restrict the parametrization of the surface. In a recent work we proposed a \emph{gauge-fixing} approach \cite{dharmavaram2021gauge, dharmavaram2022lagrangian} for fluid membranes that circumvents this issue without changing the constitutive law. Here, we extend this procedure in the presence of the embedded filaments.
\end{enumerate}

\subsection{Coupling the Rod with the Membrane}
\label{ssec:coupling}
The requirement \eqref{eq:naive constraint} that the rod remains embedded on the surface of the membrane implies that the shape of the membrane and that of the rod are inextricably coupled.  As per our assumption the lipid molecules are in a fluid state and can freely move around and across the rod. As shown in Fig.~\ref{fig:actin}(a), typically, the biofilaments are tethered to the membrane surface through membrane embedded proteins. We approximate this effect by assuming that the rod is embedded over the surface, so that the material points on the membrane surface can flow past the rod. 

\begin{figure}[htb!]
    \centering
    \begin{tabular}{cc}
    \includegraphics[width=2in]{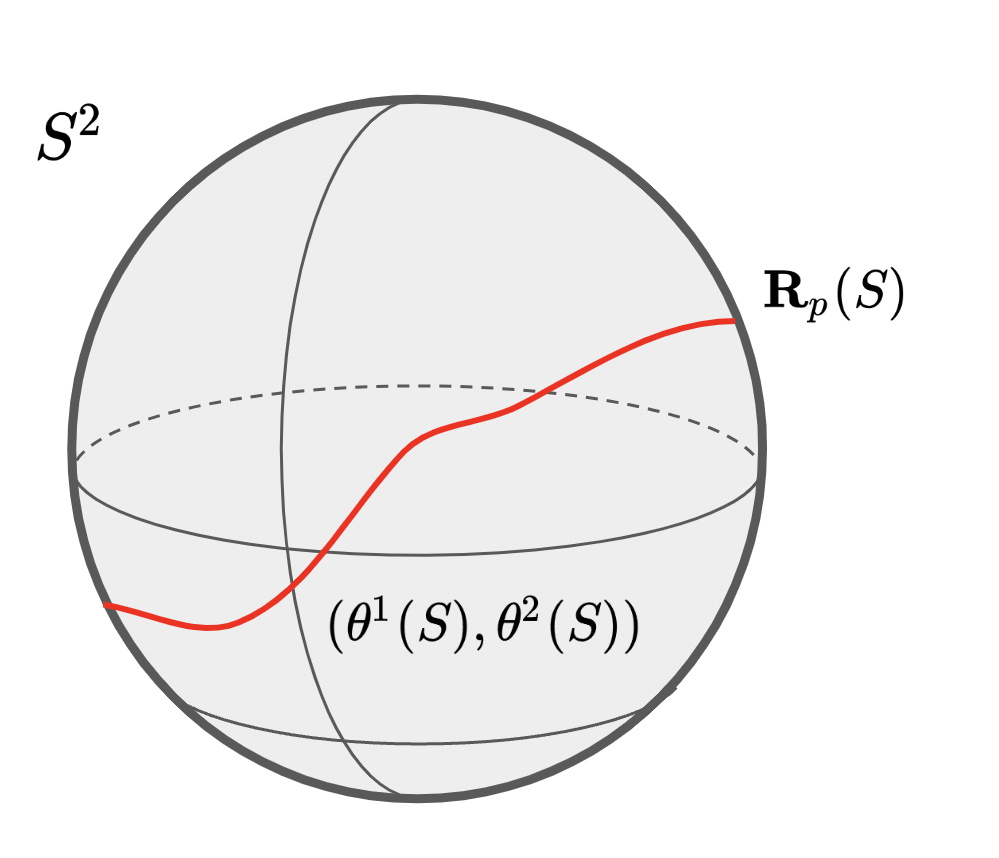} & 
    \includegraphics[width=2in]{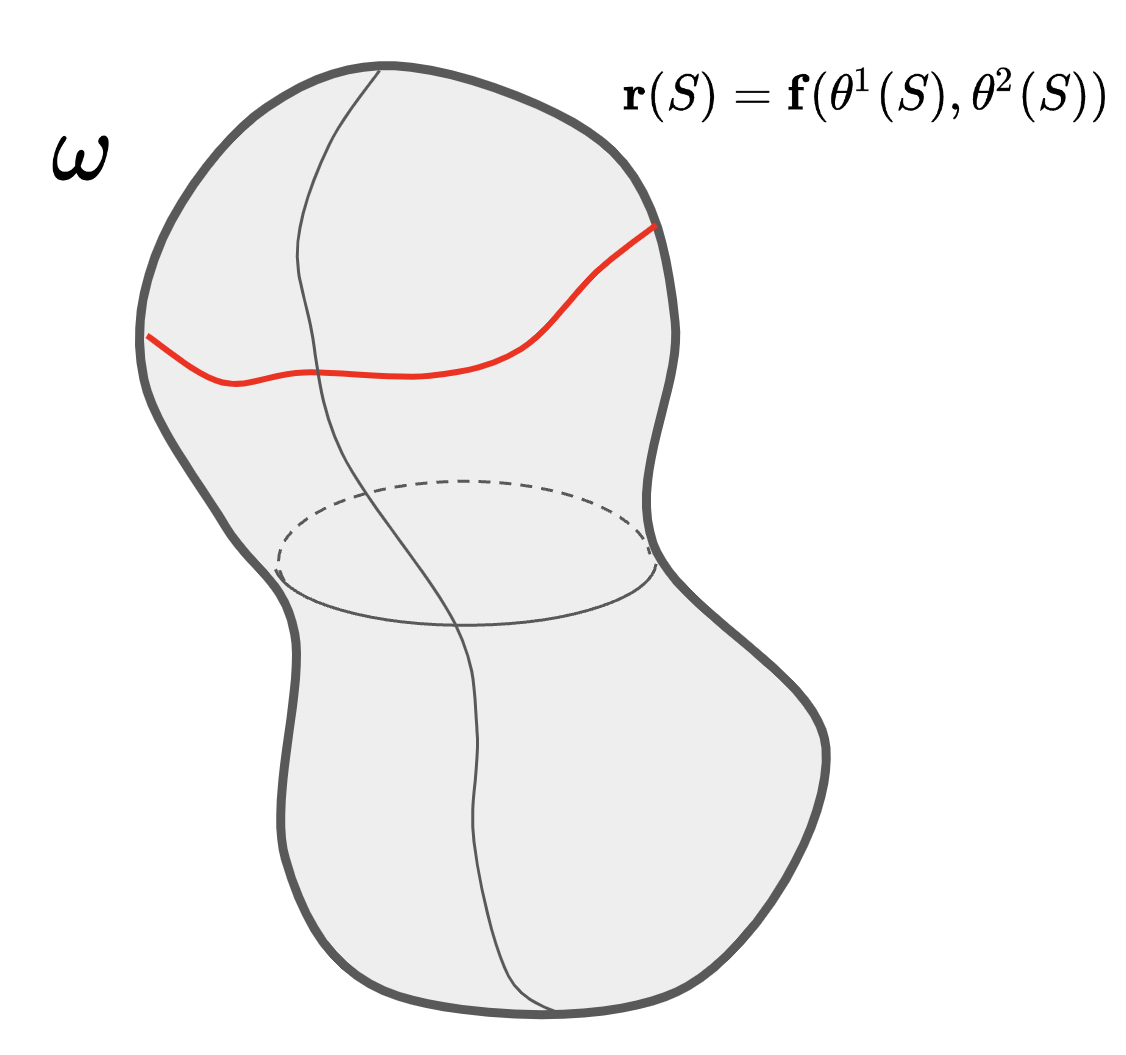}\\
    (a) & (b)
    \end{tabular}
    \caption{(a) Schematic showing a \emph{Lagrangian description} of the rod on $S^2$ in terms of spherical coordinates $\bm{\theta}(S):=(\theta^1(S),\theta^2(S))$. (b) Schematic showing the current configuration of the rod written in terms of $\mathbf{f}\circ\bm{\theta}$.  }
    \label{fig:pull back}
\end{figure}

Instead of a cumbersome constraint such as \eqref{eq:naive constraint}, we employ a \emph{Lagrangian} scheme to parametrize the deformed configuration of the rod. This idea has been inspired by our recent work on membrane-protein interaction \cite{dharmavaram2020lagrangian, dharmavaram2022lagrangian}. We introduce 
\[
\bm{\theta}:[0,L]\to {S}^2
\]
that parametrizes the deformed position of the rod onto a reference sphere $S^2$.  Instead of treating $\mathbf{r}(S)$ as the degree-of-freedom to describe the kinematics of the rod at a given point $S\in[0,L]$, we treat $\bm{\theta}(S)$ as the degree-of-freedom. For topologically spherical membranes, any choice of spherical coordinates is a natural way to parametrize $S^2$. 
Thus, the compatibility condition,
\begin{equation}
\mathbf{r}(S) = \mathbf{f}(\theta^1(S),\theta^2(S)),
\label{eq:compatibility1}
\end{equation}
automatically ensures \eqref{eq:naive constraint}, that is the rod remains embedded on the surface. Note that the right side of the previous equation can be, more abstractly, written as $(\mathbf{f}\circ\bm{\theta})(S)$. For the simplicity of notation, in above and in the rest of paper, we identify $\bm{\theta}$ with $(\theta^1,\theta^2)$, the spherical coordinates.
A schematic showing the rod on the deformed configuration and its ``pull-back'' on a reference sphere is shown in Fig.~\ref{fig:pull back}. The pull-back on the reference sphere is explicitly given by:
\[
\mathbf{R}_p(S) = \sin\theta^1(S)\cos\theta^2(S)\mathbf{i} +  \sin\theta^1(S)\sin\theta^2(S)\mathbf{j} + \cos\theta^1(S)\mathbf{k}.
\]
Note that $\mathbf{R}_p$ is \emph{not} the reference configuration of the rod. It is simply a parametrization of the deformed configuration of the rod on a fictitious unit sphere. 

While \eqref{eq:compatibility1} couples the rod's center curve to the membrane, additional conditions need to be specified to couple the rod's cross-section to the surface. In this work, this is done by mapping suitably the directors in deformed configuration to the normal vector field of the deformed surface. Without such additional conditions, the problem is ill-posed, for instance, this is so when the cross-section of rod is allowed to twist about the center curve without changing the shape of the latter. 
While many possible couplings are possible, we choose a natural scheme that is partly motivated by the biological example shown in Fig.~\ref{fig:actin}(a).  It is clear from this figure that, to a reasonable approximation, the membrane-bound protein connecting the filament to the membrane remains normal to the surface. We approximate this effect by demanding that the unit normal to the surface and a director, say $\mathbf{d}_1$, are parallel (aligned). That is,
\begin{equation}
\mathbf{d}_1(S)=\mathbf{n}(\theta^1(S),\theta^2(S)),
\label{eq:compatibility2}
\end{equation}
for all $S\in[0,L]$. 

The compatibility conditions \eqref{eq:compatibility1} and \eqref{eq:compatibility2} are sufficient to determine the other directors. Since the rod is assumed unshearable, it follows from \eqref{eq:nu_i} (for $i=3$) that
\begin{equation}
    \mathbf{d}_3(S) = \frac{\mathbf{f}_{,\alpha}(\bm{\theta})\theta^{\alpha'}(S)}{\nu_3(S)},
    \label{eq:d3=r' by nu3}
\end{equation}
for all $S\in[0,L]$, where the numerator is obtained by differentiating \eqref{eq:compatibility1}. In \eqref{eq:d3=r' by nu3}, as later in the text, we employ Einstein's summation convention over the repeated index $\alpha$. Henceforth, a summation is implied over repeated Greek indices.
Moreover, the orthonormality of directors implies
\begin{equation}
    \mathbf{d}_2(S) = \mathbf{d}_3(S)\times\mathbf{d}_1(S),
 \label{eq:d2=d3xd1}
\end{equation}
for all $S\in[0,L]$.
As a consequence of this and the compatibility equations \eqref{eq:compatibility2} and \eqref{eq:d3=r' by nu3}, $\mathbf{Q}$ is no longer a kinematic degree-of-freedom of the rod which is independent of the membrane. 

As a result of the above assumptions coupling the rod to membrane, \eqref{eq:energy 1} can be reduced to the form:
\begin{equation}
    \mathcal{E}\fnl{\mathbf{f},\bm{\theta}} =\mathcal{E}_{m}\fnl{\mathbf{f}}+\mathcal{E}_{rod}\fnl{\mathbf{f}\circ\bm{\theta}},
    \label{eq:energy 2}
\end{equation}
subject to the constraint \eqref{eq:constr}. Note that due to \eqref{eq:compatibility2}, \eqref{eq:d3=r' by nu3}, and \eqref{eq:d2=d3xd1}, the unshearability constraints \eqref{eq:unshearable} are automatically satisfied.

\subsection{Redundancy in Equilibrium Equations and Gauge Fixing}
\label{ssec:redundancy}
We now show that as a consequence of in-plane fluidity of the membrane, the Euler-Lagrange equations, i.e., first-variation of the energy, is underdetermined. In particular, we show that if the equilibrium equations of the rod hold, then the components of the equilibrium equation tangential to the membrane automatically vanish. 

\begin{figure}[htb!]
   \centering
    \begin{tabular}{c c}
       \includegraphics[width=1.5in]{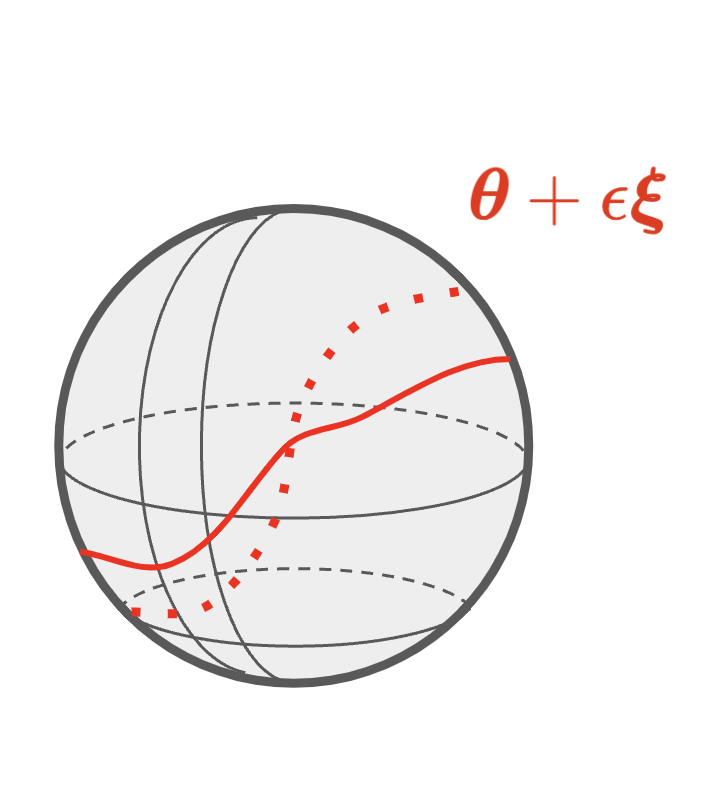}   &   \includegraphics[width=2in]{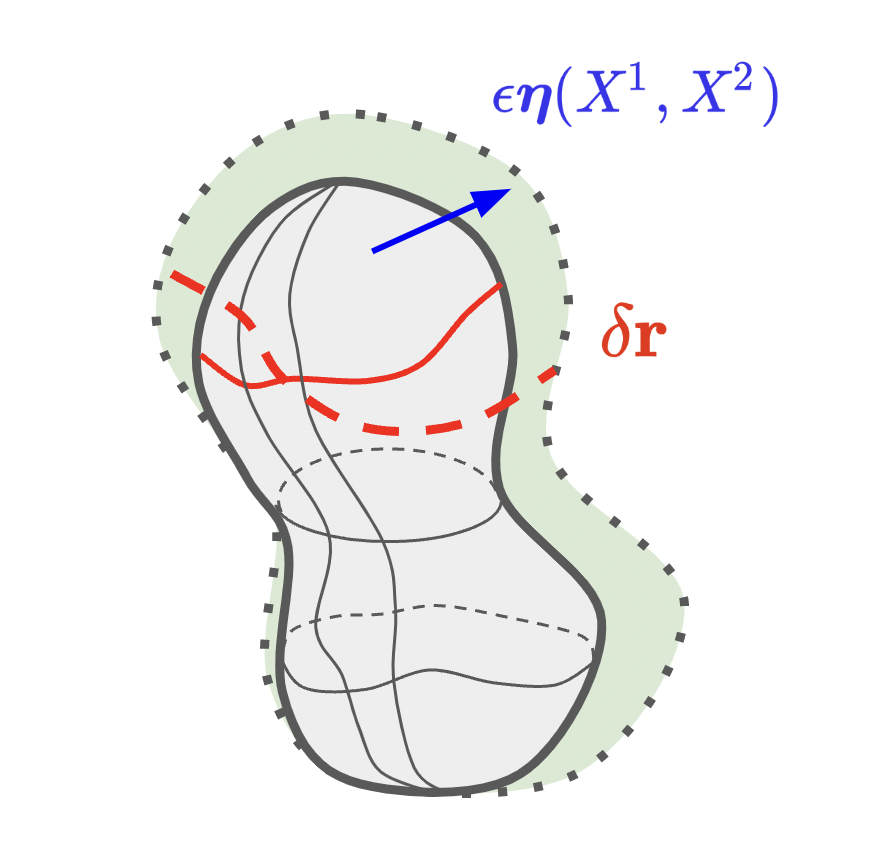}\\
       (a) & (b)
    \end{tabular}
    \caption{(a) Solid curve on the sphere represents the Lagrangian description of the rod and the dotted curve is the perturbed configuration. (b) Perturbed configuration of the deformed surface is shown in the shaded region (enclosed in the dotted curve). The perturbed current configuration of the rod is shown in the dashed curve.}
    \label{fig:perturbed configs}
\end{figure}

Following the convention of the calculus of variation, for a small parameter $\epsilon$, let us perturb the degrees of freedom as follows:
\begin{equation}
\mathbf{f}(\mathbf{X})\mapsto \mathbf{f}+\epsilon(\mathbf{v}(\mathbf{X})+w(\mathbf{X})\mathbf{n}(\mathbf{X})),
\label{eq:perturb surf}
\end{equation}
where $\mathbf{v}:=v^\alpha\mathbf{f}_{,\alpha}$ is a (smooth) tangential vector field to the surface $\omega$, $w\mathbf{n}$ defines a smooth normal vector field on $\omega$. The perturbation in $\f$ is denoted as $\epsilon\bm{\eta}$, where we define $\bm{\eta}:=\mathbf{v}+w\mathbf{n}$. Similarly, we perturb the degrees of freedom associated with the rod:
\begin{equation}
    \bm{\theta}(S)\mapsto\bm{\theta}(S) + \epsilon\bm{\xi}(S),
    \label{eq:perturb theta}
\end{equation}
for all $S\in[0,L]$ where $\bm{\xi}$ is a smooth vector field along the rod. Since $\bm{\theta}$ is the pull-back of the position vector of the rod and is defined on $S^2$, so is $\bm{\xi}$. Thus, 
\begin{equation}
\bm{\xi} = \xi^\alpha\mathbf{R}_{p,\alpha},
\label{eq:def xi1 xi2}
\end{equation}
for some scalar functions $\xi^{\alpha}$, $\alpha\in\{1,2\}$. The perturbed configurations of the membrane surface and the rod are shown schematically in Fig.~\ref{fig:perturbed configs}(a) and (b). 

Note that quantities such as $\mathbf{f}$, $\mathbf{v}$, $\mathbf{n}$, etc., are defined on the surface $S^2$ and are therefore functions of $\mathbf{X}$. On the other hand, the position vector of the rod's center curve depends on $S$.  At this point, recall the compatibility relation $\mathbf{r}(S)=(\mathbf{f}\circ \bm{\theta})(S)$, c.f.,~\eqref{eq:compatibility1}. For clarity and convenience, we introduce the following hat notation. For any vector field $\mathbf{q}:S^2\to \mathbb{R}^3$, we define 
\begin{equation}
    \hat{\mathbf{q}}:[0,L]\to\mathbb{R}^3\text{ as }\hat{\mathbf{q}} := \mathbf{q}\circ\bm{\theta}.
    \label{eq:hatnotation}
\end{equation}
For the gradient of $\mathbf{q}$, e.g., $\mathbf{q}_{,\alpha}$, we assume the notation that $\hat{\mathbf{q}}_{,\alpha} = \mathbf{q}_{,\alpha}\circ \bm{\theta}$; we employ a similar notation for higher order derivatives.

The two perturbations introduced above, c.f., \eqref{eq:perturb surf} and \eqref{eq:perturb theta}, induce the following perturbation in the position of the rod:
\begin{equation}
    \mathbf{r}(S)=\hat{\mathbf{f}}(S)=\mathbf{f}(\bm{\theta}(S)) \mapsto \mathbf{f}(\bm{\theta}(S)+\epsilon\bm{\xi}(S)) + \epsilon (\mathbf{v}+w\mathbf{n})\bigg|_{\bm{\theta}+\epsilon\bm{\xi}}.
\end{equation}
We expand the first term on the right hand side using Taylor's theorem and retain terms up to first order in $\epsilon$. Thus, we obtain
\begin{equation}
    \hat{\mathbf{f}} \mapsto \hat{\mathbf{f}} + \epsilon \Big( (\xi^\alpha+\hat{v}^\alpha)\hat{\mathbf{f}}_{,\alpha} + \hat{w}\hat{\mathbf{n}}\Big) + o(\epsilon),
    \label{eq:perturb hat r}
\end{equation}
where the first term in the outer parentheses above is the component representation of $(\nabla\hat{\mathbf{f}})^T[\bm{\xi}] + \hat{\mathbf{v}}$. Noting that $\mathbf{r}(S)=\hat{\mathbf{f}}(S)$, we see from \eqref{eq:perturb hat r} that the perturbation in the rod's deformation map has two contributions. One from the perturbation of the surface (terms involving $\hat{v}$ and $\hat{w}$ in \eqref{eq:perturb hat r}) and second from the perturbation in the Lagrangian description of the rod (term involving $\xi^\alpha$). For convenience, we denote the perturbation term---second term in \eqref{eq:perturb hat r}---by $\delta\mathbf{r}$ in Fig.~\ref{fig:perturbed configs}(b).

We introduce the following notation for the perturbations in the membrane and rod energies:
\begin{subequations}
\begin{eqnarray}
\delta_{\perp}\mathcal{E}_m[w] &:=& \frac{d}{d\epsilon}\mathcal{E}_m\fnl{\mathbf{f}+\epsilon w\mathbf{n}}\bigg|_{\epsilon=0},
\label{eq:delta Em perp}\\
\delta_{\perp}\mathcal{E}_{rod}[w] &:=& \frac{d}{d\epsilon}\mathcal{E}_{rod}\fnl{\hat{\mathbf{f}}+\epsilon \hat{w}\hat{\mathbf{n}}}\bigg|_{\epsilon=0},
\label{eq:delta Er perp}\\
\delta_{\parallel}\mathcal{E}_m[\mathbf{v}] &:=& \frac{d}{d\epsilon}\mathcal{E}_m\fnl{\mathbf{f}+\epsilon\mathbf{v}}\bigg|_{\epsilon=0},
\label{eq:delta Em ||}\\
\delta_{\parallel}\mathcal{E}_{rod}[q^1,q^2] &:=& \frac{d}{d\epsilon}\mathcal{E}_{rod}\fnl{\hat{\mathbf{f}}+\epsilon (q^1\hat{\mathbf{f}}_{,1}+q^2\hat{\mathbf{f}}_{,2})}\bigg|_{\epsilon=0},
\label{eq:delta Er ||}
\end{eqnarray}
\label{eqs:delta Es}\end{subequations}
where $\delta_{\parallel}\mathcal{E}_{m}$, etc. are linear operators of their arguments (indicated in brackets). Note that $\delta_\perp$ and $\delta_{||}$ denote variations due to perturbations normal and tangential (respectively) to the membrane.

To derive the Euler-Lagrange equations, we set the first variation of \eqref{eq:energy 2}, due to the perturbations \eqref{eq:perturb surf} and \eqref{eq:perturb hat r}, to zero. That is,
\begin{equation}
\frac{d}{d\epsilon}\mathcal{E}\bigg|_{\epsilon=0} = (\delta_\perp \mathcal{E}_m+\delta_\perp \mathcal{E}_{rod})[w] + (\delta_{\parallel} \mathcal{E}_m+\delta_{\parallel} \mathcal{E}_{rod})[v^1,v^2] + \delta_{\parallel} \mathcal{E}_{rod}[\xi^1,\xi^2] = 0.
\end{equation}

Since the perturbations (enclosed in square brackets) are independent (and arbitrary), we obtain the following:
\begin{subequations}
\begin{eqnarray}
    (\delta_\perp \mathcal{E}_m+\delta_\perp \mathcal{E}_{rod})[w] &=& 0,
    \label{eq:first var norm}\\
    (\delta_{\parallel} \mathcal{E}_m+\delta_{\parallel} \mathcal{E}_{rod})[v^1,v^2] &=& 0,
     \label{eq:first var tgt}\\
    \delta_{\parallel} \mathcal{E}_{rod}[\xi^1,\xi^2] &=& 0.
     \label{eq:first var rod}
\end{eqnarray}
\label{eqs:first vars}
\end{subequations}
Here, \eqref{eq:first var norm} can be interpreted as the normal component of the momemtum balance for the membrane, \eqref{eq:first var tgt} as the tangential component of the momentum balance, and \eqref{eq:first var rod} as the rod's momentum balance.

It is a well-known fact \cite{capovilla2003deformations, steigmann2003variational, dharmavaram2021gauge} that the tangential variation of the Helfrich-Canham energy is identically zero. That is, $\delta_{\parallel}\mathcal{E}_m[\mathbf{v}]\equiv 0$, for all tangential perturbations $\mathbf{v}$, and on any configuration defined by the deformation map $\mathbf{f}$. Thus, if \eqref{eq:first var rod} holds, then \eqref{eq:first var tgt} is trivially zero. Thus, we do not obtain any new condition for the tangential equilibrium equation for the membrane. In other words, the equilibrium condition for the membrane is underdetermined.

\subsubsection{The Gauge-fixing Procedure}
\label{sssec:gauge}
We now address the second computational challenge mentioned at the beginning of this section. A strategy to circumvent problems arising from membrane's in-plane fluidity has been presented in recent works \cite{dharmavaram2021gauge, dharmavaram2022lagrangian} through a procedure termed \emph{gauge-fixing}. While these earlier work considered only a fluid membrane, here we demonstrate that the procedure can be employed for the membrane-rod system. Note that the crucial ingredients for this argument are that (1) the surface has a spherical topology, and (2) the rod floats on the membrane. We emphasize that the gauge-fixing procedure is \emph{not} the main contribution of the paper. We discuss it simply to convince the reader that the inclusion of a rod in the system, in the manner described above, does not preclude the gauge-fixing procedure that was developed solely for fluid membranes.

As we noted above, the equilibrium equations for the membrane are underdetermined. The main idea of \emph{gauge-fixing} is to supplement the equilibrium equations of the membrane with another constraining equation. A classical example of this procedure is in potential-theoretic formulation of electrodynamics. Since vector and scalar potentials cannot be uniquely determined, extra-conditions called \emph{gauge conditions} are imposed on the potentials (e.g., Coulomb gauge that sets the divergence of the vector potential to be zero). For the Helfrich model on topologically spherical surfaces, a recently developed procedure has been presented in   \cite{dharmavaram2021gauge}. A detailed justification of this procedure has been provided in this reference. Below, we only discuss the procedure from an application viewpoint. 

We supplement the total energy of the system with the \emph{harmonic map energy}, defined by
\begin{equation}
    \E_{\text{HM}}\fnl{\mathbf{f}} = \int_{\omega} \frac{1}{2}{g}^{\alpha\beta}h_{\alpha\beta}\;da,
    \label{eq:har en conf cor}
\end{equation}
where $g_{\alpha\beta}$ ($\alpha,\beta=1,2$) are the (contravariant) components of the metric tensor of the deformed surface $\omega$, and $h_{\alpha\beta}$ are the (covariant) components of the metric tensor of in a chosen parametrization of the reference surface $S^2$. Note that $h_{\alpha\beta}$ is fixed for a given parametrization of $S^2$ (see \cite{dharmavaram2021gauge} for an example involving spherical coordinates). On the other hand $g_{\alpha\beta}$ can be computed in terms of the deformation map via:
\begin{equation}
g_{\alpha\beta}=\mathbf{f}_{,\alpha}\cdot\mathbf{f}_{,\beta},\;\alpha,\beta\in\{1,2\},
\end{equation}
where the commas represent partial derivatives with respect to $X^1$ and $X^2$. Thus, the required equilibria are computed as critical points of the new functional:
\begin{equation}
    \tilde{\mathcal{E}}\fnl{\mathbf{f},\bm{\theta}}= \mathcal{E}_m\fnl{\mathbf{f}}+\mathcal{E}_{rod}\fnl{\mathbf{f}\circ\bm{\theta}}+\mathcal{E}_{\text{HM}}\fnl{\mathbf{f}}.
    \label{eq:total energy}
\end{equation}
We term this the \emph{gauge-fixed formulation}.

A nontrivial feature of the gauge-fixing procedure is that adding the harmonic energy \eqref{eq:har en conf cor} does not alter the constitutive law of the fluid membrane. This is despite the fact that $\mathcal{E}_{\text{HM}}$ depends on the deformation map of the surface. It can be shown \cite{dharmavaram2021gauge, dharmavaram2022lagrangian} that the Euler-Lagrange operator corresponding to $\mathcal{E}_{\text{HM}}$ only contributes in the direction tangential to the surface (see Theorem 3, Appendix D of \cite{dharmavaram2021gauge}). The first variation of \eqref{eq:total energy} therefore supplements the already established conditions \eqref{eqs:first vars} with the additional condition $\delta_{||}\mathcal{E}_{\text{HM}}[v^1,v^2]=0$.  Thus, we do not modify the physics of the problem event though we supplement the total energy with \eqref{eq:har en conf cor}. Other established methods~\cite{feng2006finite, ma2008viscous, Elliot} modify the constitutive law by adding additional terms to the energy. Typically, viscous terms that penalize in-plane motions are added to the energy and the equilibrium configurations are computed by iteratively decreasing the effect of the added terms. Since the gauge-fixed formulation does not require this iterative step, our computations are significantly faster.

One notable feature of the gauge-fixing procedure presented above is that adding $\E_{\text{HM}}$ is not sufficient to constrain all the tangential degrees of freedom. That is, the tangential equilibrium equations are still underdetermined. We have shown in \cite{dharmavaram2021gauge} that the null space of the linearized Harmonic Map equation (i.e., $\delta_{||}\mathcal{E}_{\text{HM}}[v^1,v^2]=0$) is six-dimensional and the null vectors can be generated by M\"obius group of transformations. The six linearly independent modes can be visualized as the rigid translational and rotational modes of the sphere~\cite{Arnold_Rogness}. We fix the translational modes via the ``zero-mass'' constraint~\cite{dharmavaram2021gauge}:
\begin{equation}
  \int_{\omega} \mathbf{R}(X^1,X^2) \;da = \mathbf{0}\, ,
  \label{eq:zero mass}
\end{equation}
where $\mathbf{R}(X^1,X^2)$ is a parametrization of the reference sphere. While the three rotational modes could be fixed by imposing additional 
``landmark constraints''~\cite{Gu_Wang}, for the results presented below, we found that this was not necessary. The minimization 
procedure
was able to find reliable solutions without any additional constraints.

\section{First Variation and Weak Form}
\label{sec:weak form}
In this section we derive the weak-form for the gauge-fixed formulation, cf.~\eqref{eq:total energy}, subjected to constraint \eqref{eq:constr} and  \eqref{eq:zero mass}. We begin by setting the notation and recalling basic facts about differential geometry that we use to derive the weak form. 

The tangent space of $\omega$ at $\mathbf{X}$ (parametrized by $X^1$ and $X^2$) is spanned by the basis 
\begin{equation}
  \mathbf{a}_\alpha := \f_{,\alpha} = \frac{\partial\f}{\partial X^\alpha}\in\mathbb{R}^3,\;\alpha\in\{1,2\}.
\end{equation}
The components of the metric tensor of $\omega$ are given by
$$g_{\alpha\beta}=\mathbf{a}_\alpha\cdot\mathbf{a}_\beta,$$
where dot represents the Euclidean inner product of vectors in $\mathbb{R}^3$. We denote the dual basis vectors superscripts, i.e.,
$$\mathbf{a}^\alpha := g^{\alpha\beta}\mathbf{a}_\beta,$$
where $g^{\alpha\beta}$ are the contravariant components of the metric tensor. It follows that the unit normal field to $\omega$ is given by 
$$\mathbf{n} = \frac{\mathbf{a}_1\times \mathbf{a}_2}{\sqrt{g}},$$ 
where we have used the identity $\sqrt{g}=\parallel{}\mathbf{a}_1\times\mathbf{a}_2{}\parallel$. The components of the second fundamental form of $\omega$ are
$$b_{\alpha\beta}=-\mathbf{n}_{,\alpha}\cdot \mathbf{a}_\beta = \mathbf{n}\cdot\mathbf{a}_{\alpha,\beta}.$$

The Mean and Gaussian curvatures are given by
\begin{equation}
  H=\frac{1}{2}b^\alpha_\alpha =
  -\frac{1}{2}\mathbf{a}^\alpha\cdot\mathbf{n}_{,\alpha}\, , \quad   K=\det(b^\alpha_\beta)\, ,
  \label{eq:HK}
\end{equation}
where, following the standard convention, the tensor indices are lowered or raised by multiplying using the metric tensor $g_{\alpha\beta}$ or its inverse $g^{\alpha\beta}$, respectively. 

To derive the weak form of the gauge-fixed functional \eqref{eq:total energy}, we compute its variation $\delta\tilde{\E}$ with respect to the variations \eqref{eq:perturb surf} and \eqref{eq:perturb theta}. Equilibrium configurations are obtained by the first-variation condition:
\begin{equation}
    \delta\tilde{\E} = \delta\E + \delta\E_{\text{HM}} + \delta\E_{rod}=0,
    \label{eq:weak gauge fixed}
\end{equation}
subject to constraints \eqref{eq:constr} and \eqref{eq:zero mass}.
It can be shown (see~\cite{ma2008viscous, dharmavaram2021gauge, dharmavaram2022lagrangian} for details) that:
\begin{equation}
  \delta\E_m = \int_{S^2}
  \left(\mathbf{n}^\alpha\cdot\delta\mathbf{a}_\alpha +
    \mathbf{m}^\alpha\cdot\delta\mathbf{n}_{,\alpha}-\frac{1}{3}pV\mathbf{n}\cdot\bm{\eta}\right)\;JdA,
  \label{eq:weak_helf}
\end{equation}
\begin{equation}
    \delta\E_{\text{HM}} = \frac{1}{2}\int_{S^2}h_{\alpha\beta}\Big(-2g^{\alpha\gamma}
\mathbf{a}^\beta+g^{\alpha\beta}\mathbf{a}^\gamma\Big)
\cdot\bm{\eta}_{,\gamma} J\;dA,
\label{eq:weak harmonic energy}
\end{equation}
where $J=\sqrt{g}/\sqrt{h}$ with $h$ being the determinant of $h_{\alpha\beta}$, $\bm{\eta}$ is the variation in the surface (see note following \eqref{eq:perturb surf}), and 
\begin{multline}
  \mathbf{n}^\alpha = (\kappa
  (H-C_0)+2\kappa_g H)g^{\alpha\beta}\mathbf{n}_{,\beta}+(\kappa(H-C_0)^2+\kappa_g(K-b^\nu_\beta \mathbf{a}^\beta\cdot\mathbf{n}_{,\nu}))\mathbf{a}^\alpha\\
  -\frac{1}{3}pV((\mathbf{f}\cdot\mathbf{n})\mathbf{a}^\alpha - (\mathbf{f}\cdot\mathbf{a}^\alpha)\mathbf{n})\, ,
  \label{eq:n}
\end{multline}
\begin{equation}
  \mathbf{m}^\alpha = \left(-\kappa(H-C_0)\delta^\alpha_\beta+\kappa_g(b^\alpha_\beta - 2H\delta^\alpha_\beta)\right)\mathbf{a}^\beta\, ,
  \label{eq:m}
\end{equation}
where $\delta^{\alpha}_\beta$ is the Kronecker delta, and the variations appearing in the above equations are given by \cite{feng2006finite}
\begin{equation}
  \delta\mathbf{a}_\alpha = \bm{\eta}_{,\alpha},
  \label{eq:delta_a}
\end{equation}
\begin{equation}
  \delta\mathbf{n}_{,\alpha} = -(\mathbf{a}^\beta_{,\alpha}\otimes \mathbf{n} + \mathbf{a}^\beta\otimes\mathbf{n}_{,\alpha})\cdot\bm{\eta}_{,\beta} - (\mathbf{a}^\beta\otimes\mathbf{n})\cdot\bm{\eta}_{,\beta\alpha}.
  \label{eq:delta_d_alpha}
\end{equation}

The first variation of $\E_{rod}$ is given by
\begin{equation}
    \delta\E_{rod} = \int_0^L D(\nu_3-1)\delta \nu_3 + E \sum_{i=1}^3 \kappa_i\delta\kappa_i\;dS,
    \label{eq:weak rod}
\end{equation}
where the explicit formulas for $\delta\nu_3$, $\delta\kappa_i$ ($i=1,2,3$) in terms of $\bm{\eta}$ and $\bm{\xi}$ are provided in the \ref{sec:app:deriv fv strains}; see Propositions \ref{prop:nu3}-\ref{prop:k3}.
Anticipating the key numerical computations, we conclude this section with the following remark. 
\begin{remark}
It has been noted in \cite{dharmavaram2022lagrangian} that computing equilibria of  $\tilde{\E}$ using its weak form \eqref{eq:weak gauge fixed} is not always feasible with a gradient-descent based minimization solver. This is because $\E_{\text{HM}}$ does not need to be positive, and so a minimizer of $\tilde{\E}$ need not be that of $\E$. Thus, the gauge-fixing procedure potentially could transform local minima into saddle points and local maxima. Note that they continue to remain critical points. To deal with this issue in computations, we instead minimize the following functional:
\begin{equation}
    \hat{\E}:=\mathcal{E}_m + \mathcal{E}_{rod}+\lambda_g(\mathcal{E}_{\text{HM}})^2 + \mu_1 \Big(\int_{S^2}J\;dA-4\pi\Big)^2 + \mu_2\Big|\int_{S^2} \mathbf{R}(X^1,X^2) J\;dA\Big|^2\,,
    \label{eq:obj fun}
\end{equation}
where $\lambda_g$ can be interpreted as a penalty parameter that controls the strength of the harmonic map energy. As can be seen in the last two terms of \eqref{eq:obj fun}, we enforce the constraints \eqref{eq:constr} and \eqref{eq:zero mass} through penalty constraints with parameters $\mu_1$ and $\mu_2$. The associated weak form is given by:
\begin{equation}
    \delta\E_m + \delta\E_{rod} + (2\lambda_g\E_{\text{HM}}) \delta\E_{\text{HM}} + \lambda_1 \int_{S^2} \mathbf{a}^\alpha\cdot\bm{\eta}_\alpha J\;dA + \lambda_2\int_{S^2}\mathbf{R}(X^1,X^2) \mathbf{a}^\alpha\cdot\bm{\eta}_\alpha J\;dA= 0 \, ,
    \label{eq:weak form modified functional}
\end{equation}
where $\lambda_1=2\mu_1(\int_{S^2}J\;dA-4\pi)$ and $\lambda_2=2\mu_2 \int_{S^2}\mathbf{R}J\;dA$. Note that we have used $\delta J = J\mathbf{a}^\alpha\cdot\bm{\eta}_\alpha$ to obtain the last two terms of \eqref{eq:weak form modified functional} while the variations $\delta\E_m$, $\delta\E_{\text{rod}}$, and $\delta\E_{\text{HM}}$ are given by \eqref{eq:weak_helf}, \eqref{eq:weak rod}, and \eqref{eq:weak harmonic energy}, respectively.
\end{remark}

\section{Numerical Results}
\label{sec:numerical results}
\subsection{Finite Element Discretization}
\label{ssec:loop fem}

We employ a Ritz-Galerkin projection scheme to discretize the energy \eqref{eq:obj fun} and its weak form \eqref{eq:weak form modified functional}. Since our system is a combination of two distinct continuua---a shell and a rod---we use two different finite elements to discretize them. We use a Loop subdivision finite element \cite{cirak2000subdivision} for the fluid shell and a Hermite cubic element for the rod. 

The Loop subdivision finite element approach has been effectively applied to solid and fluid shells \cite{cirak2000subdivision, feng2006finite, ma2008viscous, dharmavaram2021gauge, dharmavaram2022lagrangian}. While there are a variety of finite elements available for thin shells, we use the Loop subdivision method scheme as it is relatively straightforward and only uses the positional degrees-of-freedom for the mesh. In this method, the surface is discretized into triangles, and the deformation map $\mathbf{f}$ is interpolated on each triangular element using bivariate quartic box-splines, expressed in terms of the barycentric coordinates of the element. The method guarantees $C^1$ continuity across elements which is required for discretizing $\mathcal{E}$. The trade-off with using only the nodal degrees of freedom is that unlike traditional finite elements that have test/shape functions supported only on the element, in the Loop shell scheme, the functions (box-splines) are supported on the first ring of neighbors around a given element. Explicit implementation details on this method can be found in \cite{cirak2000subdivision, feng2006finite, ma2008viscous}. In what follows, we denote $N_{\text{surface}}$ as the number of nodes of the mesh. Each node contains three degrees-of-freedom for the $x$-, $y$-, and $z$-components of the deformation map.

The relevant degrees of freedom for the rod are $\theta^1(S)$ and $\theta^2(S)$. The anchoring conditions for the rod \eqref{eq:compatibility1} and \eqref{eq:compatibility2} means that the directors of the rod are automatically determined from the deformation map $\mathbf{f}$, $\theta^1(S)$, and $\theta^2(S)$. Since the bending and twisting energy of the rod depends on the second derivative of $\mathbf{r}$, and hence that of $\theta^\alpha(S)$ ($\alpha\in\{1,2\}$), we use Hermite cubic elements to discretize the rod. We discretize the reference arc length of the rod $S\in[0,L]$ into $N_{\text{rod}}$ nodes. For each node and for each $\alpha\in\{1,2\}$, we have two degrees-of-freedom corresponding to the Hermite cubic elements. The rod's energy ($\mathcal{E}_{rod}$) and first variation ($\delta\mathcal{E}_{rod}$) depend on $\hat{\mathbf{f}}(S)$, cf., \eqref{eq:hatnotation}. We compute the discretized integrals using a three-point Gaussian quadrature. To evaluate $\hat{\mathbf{f}}$ at the quadrature points, we have to compute the shape functions of the Loop subdivision element at $(\theta^1,\theta^2)$. We first use a search algorithm developed in our earlier work \cite{dharmavaram2022lagrangian} to determine the corresponding element of the surface mesh and then use quasi-Newton solver to determine the corresponding barycentric coordinates of the quadrature point in the element. 

We present computational results obtained using the discretization discussed above. We use a gradient-based L-BFGS minimization routine \cite{zhu1997algorithm} to minimize \eqref{eq:obj fun}. Unless otherwise stated, the minimization routine was called with a randomized initial configuration. We set the penalty parameters (cf.,~\eqref{eq:obj fun}) to be $\mu_1=\mu_2=10^3$ and the preferred curvature $C_0=0$. We limit our work to spherical topologies. As a consequence of the Gauss-Bonnet theorem, the term involving  $c_g$ in \eqref{eq:energy_pre} is a constant and thus does not affect the equilibrium configuration. The parameter $\lambda_g$ appears as a penalty parameter for the harmonic map energy in \eqref{eq:obj fun}. A small value for $\lambda_g$ can cause mesh distortion (similar to the ones noted in \cite{feng2006finite}), while a large value can slow down convergence. In the results presented below, we choose its value (by experimentation) so as to aid convergence. A value of $\lambda_g=50$ worked for most cases.

\subsection{Verification and Validation}
\label{ssec:validation}
In this section, we validate our computational formulation with two studies: a convergence test and a benchmark study.

\begin{figure}[htb!]
    \centering
    \includegraphics[width=0.95\textwidth]{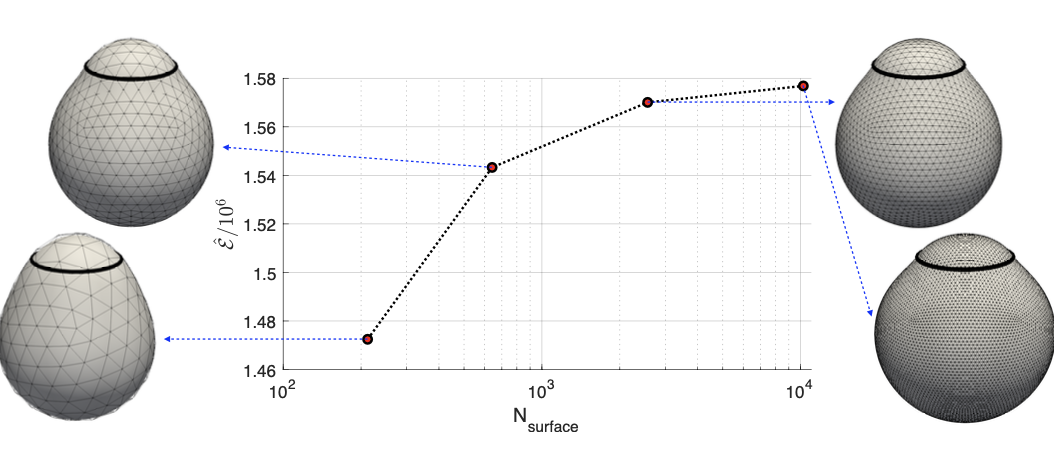}
    \caption{Convergence of discretized energy $\hat{\mathcal{E}}$ versus number of nodes in the Loop surface mesh $N_{\text{surface}}$. The mesh in the equilibrium deformed configuration is also shown corresponding to the selected cases.}
    \label{fig:convergence}
\end{figure}

Fig.~\ref{fig:convergence} summarizes our results from the convergence study. In the central plot, we have presented the equilibrium state's total (discretized) energy $\hat{\mathcal{E}}$ (normalized by a factor of $10^6)$ on the y-axis versus (logarithm of) $N_{\text{surface}}$ on the x-axis, for four different discretizations of the surface mesh. We used $N_{\text{surface}}=212,\;642,\;2562,\;10242$ with $D=100$, $E=1$, $p=10$, $N_{rod}=31$, and $c=5$. The reference configuration for the rod was chosen to be a circular state with $L=\pi$ (and $\alpha_{31}=-1,\;\alpha_{32}=0$; see \eqref{eq:D1 D2 D3}). Note the convergence in the energy as we increase the degrees-of-freedom. The corresponding equilibrium solutions are shown as insets.

In our second validation study, we compare the results obtained using our computational formulation with a semi-analytical benchmark problem of ``line tension'' induced deformation in the membrane as considered by Das et al in \cite{das2008higher}. In their work, the authors considered spherical lipid membranes, parts of which can exist in two different fluid phases. The boundary between the two phases (which is a closed curve) creates a line tension which causes the formation of a neck. In Sec. 5.1 of \cite{das2008higher}, the authors provide expressions for a special case, where the two phases have the same material properties and the phase boundary lies on the equatorial plane. Here, we approximate this result by modeling the filament as a flexible string with a large elastic modulus $D=100$ (incapable of supporting bending, i.e., $E=0$) compared to the bending modulus of the membrane $c=1$. In this scenario the mechanics of rod is essentially that of a rubber band, thus approximating the effect of line tension. The filament is further constrained to lie on the equatorial plane by setting $\theta^1(S)=\pi/2$ in the L-BFGS solver. The membrane is pressurized from the inside with a relatively large pressure $p = 1000$ (compared to $c=1$).  In \cite{das2008higher}, a perturbation expansion for the solution in terms of the small parameter $\epsilon:=1/\sqrt{p}$ has be obtained. More specifically, the authors provide the following relationship between the line tension in the string ($\sigma$), as a function of the neck radius ($r^*$):
\begin{equation}
    \sigma = \sigma_0 + \epsilon \sigma_1\, ,
    \label{eq:sigma}
\end{equation}
where 
\begin{equation}
    \sigma_0 = \frac{(r^* - r^{*3})}{\sqrt{2 - r^{*2}}},\;\sigma_1 = ( -0.98517 + 3.36358 r^* - 2.15325 r^{*2} + 0.84090 r^{*3} - 0.93567 r^{*4})\, ,
    \label{eq:sigma0 sigma1}
\end{equation}
and $r^*$ is the neck radius. We compute the tension in the filament using the following equation\begin{equation}
    \sigma=\frac{D}{p}(\bar{\nu}_3-1)\, ,
\end{equation}
where $\bar{\nu}_3$ is the average $\nu_3$ in along the filament.

\begin{figure}[htb!]
    \centering
    \includegraphics[width=0.95\textwidth]{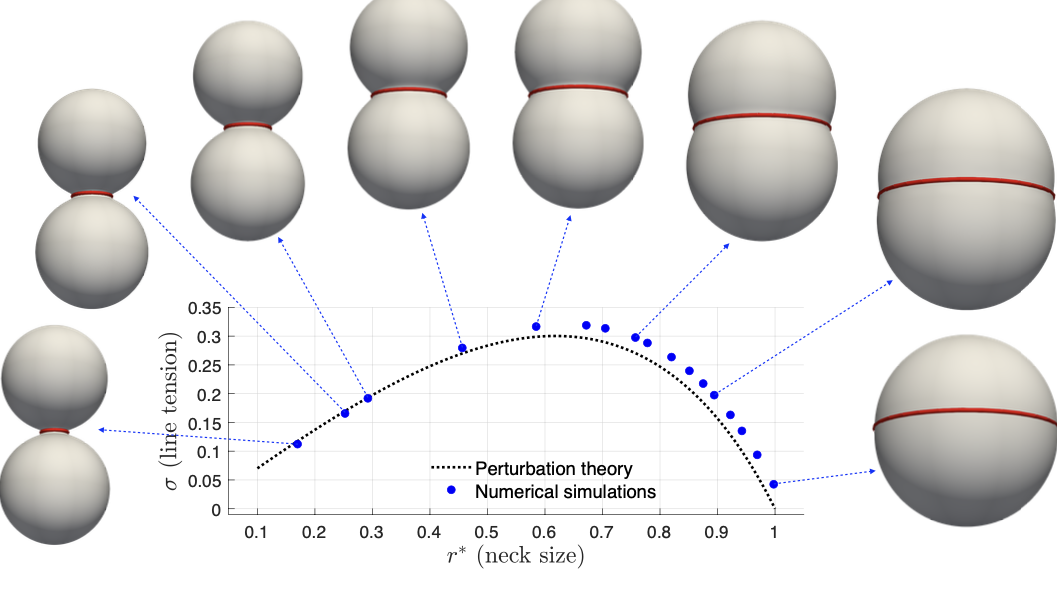}
    \caption{Numerical simulations versus perturbation theory solutions as a function of $r^*$. From left to right the numerical simulations corresponding to the circular markers are obtained for $L/2\pi = \{0.08, 0.095, 0.10, 0.12, 0.14, 0.16, 0.17, 0.19, 0.2, 0.225, 0.25, 0.275, 0.3, 0.35, 0.4, 0.5, 0.7\}$. Inserts are not necessarily to scale.     }
    \label{fig:validation}
\end{figure}

Comparison between the numerical results and the semi-analytical results is shown in Fig.~\ref{fig:validation}. The dotted curve is the perturbation solution for $\sigma$ as a function of the neck radius $r^*$ computed with \eqref{eq:sigma} and \eqref{eq:sigma0 sigma1}. The scatter circular markers correspond to the solutions computed using the proposed finite element formulation for various reference configurations of the rod, which we set to be a circle with length between $L/2\pi=0.08$ to $0.7$. Note that due to high internal pressure and  extensiblity of the rod, we find that $L/2\pi\neq r^*$.

\subsection{Applications}
We now apply the framework, developed thus far, to three applications. 

\begin{figure}[htb!]
     \centering
    \includegraphics[width=0.95\textwidth]{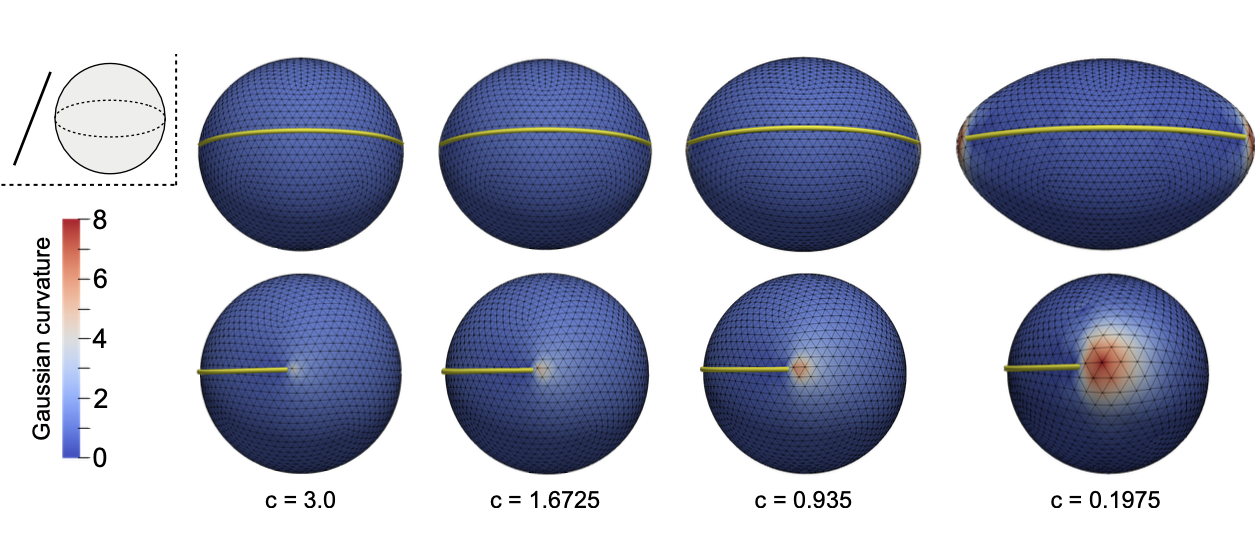}
    \caption{Equilibrium configurations of a straight (in the reference configuration) rod with length $L=\pi$ on a spherical membrane as a function of membrane bending stiffness $c$ from high (left) to low (right). The deformed rod is represented in yellow while the membrane surface is colored by its Gaussian curvature.}
    \label{fig:straight rod on sphere}
\end{figure}

In the first application, we consider a straight rod, i.e., whose reference configuration is flat, embedded on the surface of a spherical membrane. We represent this schematically in the inset at the top left corner of Fig.~\ref{fig:straight rod on sphere}. For a straight rod, it follows from \eqref{eq:D'} that $\alpha_{31}=\alpha_{32}=0$. The same figure shows the equilibrium configurations when length $L=\pi$ as we vary the membrane stiffness ($c$). Observe that when $c$ is large, it is more energetically favorable for the rod to undergo bending as it remains embedded on the membrane and the configuration is close to being spherical. As $c$ is decreased, the rod straightens out and the configuration becomes progressively more bean-shaped. For these simulations, we set $D=20$, $E=2$, and $p=0$.

\begin{figure}[htb!]
\centering
    \includegraphics[width=0.95\textwidth]{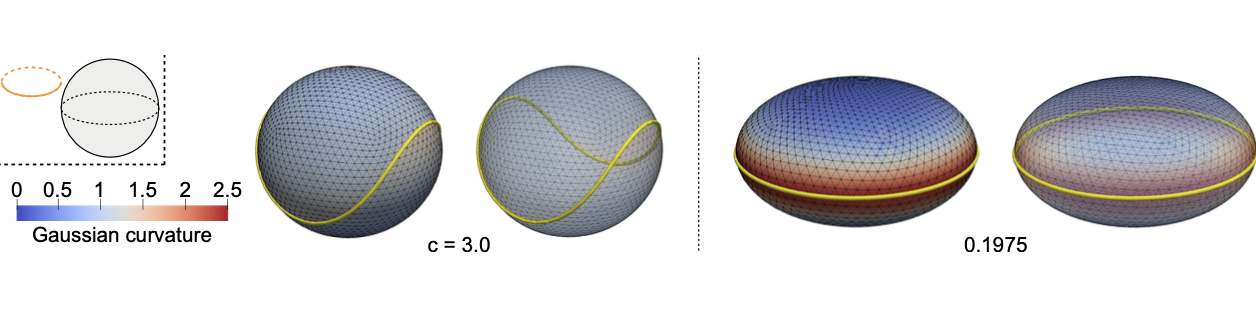}
    \caption{Coupling between membrane and rod bending. If the bending stiffness of the membrane is large, the rod buckles (achieving a configuration that recalls the seam of a baseball) while the membrane configuration remains close to spherical (left, $c=3.0$). As the membrane bending stiffness decreases, the rod maintains its circular (close to reference) configuration while the membrane tends toward a ``pancake'' configuration (right, $c=0.1975$). In both cases, the membrane deformed configurations are rendered with an opaque and semi-transparent surface to clearly show the rod deformed configuration on the entire membrane. For these simulations, we set $D=20$, $E=1$, $p=0$.}
    \label{fig:circularRodSphere}
\end{figure}

As a second application, we next consider a rod whose reference configuration is a circular ring as represented schematically in the inset in Fig.~\ref{fig:circularRodSphere}.  Note that the ring freely floats on the surface of the membrane. In this simulation, the length of the rod is set to be greater than the circumference of the unit sphere ($L=2.6\pi > 2\pi$). If the bending stiffness of the membrane is large (compared to that of the rod), the rod undergoes buckling, while the membrane continues to remain almost spherical. This configuration is shown in Fig.~\ref{fig:circularRodSphere} (left) where $c=3.0$. The rod's configuration can be visualized as the seam of a baseball, as seen in the figure. For smaller membrane stiffness (e.g, $c=0.1925$ as shown in Fig.~\ref{fig:circularRodSphere}, right), a ``pancake'' configuration is favored. In this configuration, the rod remains circular and close to its referential configuration, but the membrane undergoes large deformation. For these simulations, we have used $D=20$, $E=1$, $p=0$.

\begin{figure}[htb!]
    \centering
    \includegraphics[width=0.95\textwidth]{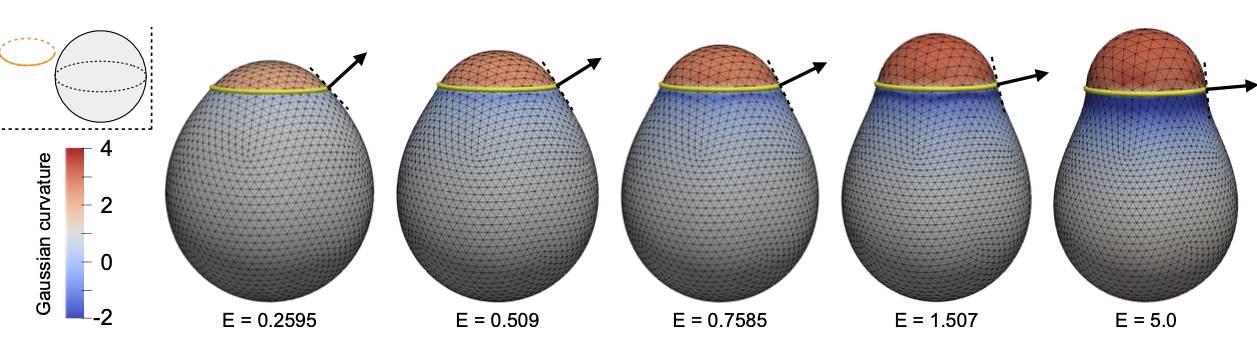}
    \caption{Coupling between membrane bending and rod twisting: as the rod stiffness $E$ increases from left to right, the rod `deforms' the membrane to maintain its director $\mathbf{d}_1$ in the horizontal reference configuration. In this example, the membrane bending stiffness $c=1$ and the rod stretching modulus $D=10$. The black arrow depicts the normal to the membrane surface at the rod contact location and hence the rod director $\mathbf{d}_1$.}
    \label{fig:bendtwist}
\end{figure}

In our third application, we explore the coupling between membrane bending and rod twisting. We consider a rod with a circular referential configuration with $L=\pi$ embedded on the membrane, see the schematic inset at the top left corner of Fig.~\ref{fig:bendtwist}. The referential director $\mathbf{D}_1$ points radially outward in the plane of the ring, cf., \eqref{eq:D1 D2 D3}(a). We fix $c=1$, $D=10$, and $p=0$. As we increase rod bending/twisting stiffness, $E$, the bud at the top end of the membrane becomes more prominent, see Fig.~\ref{fig:bendtwist}. This behavior can be understood as follows. For small values of $E$, the twisting in the rod is more predominant compared to bending in the membrane. This can be seen in the first two figures at the left end of Fig.~\ref{fig:bendtwist} where the director $\mathbf{d}_1$ of the deformed state (shown as the black arrow) is not in the plane of the ring. As $E$ increases, the bending in the membrane becomes prominent and at the same time the rod does not twist as much. Note that for large $E$, the director $\mathbf{d}_1$ is close to pointing radially outward in the plane of the ring. We used $N_{\text{rod}}=31$ and $N_{\text{surface}}=2562$ for these simulations.

\section{Discussion and Conclusions}
\label{sec:discussion and conclusions}
In this work, we propose a computational strategy to understand the coupled mechanics of a semiflexible filament embedded, yet freely floating, on a fluid membrane. We model the membrane and the filament as continua and use a variational approach to compute equilibrium configurations of the system. The requirement that the filament remains embedded poses computational challenges that we address in this work. We do this by using a Lagrangian approach where the rod is parametrized by functions $\theta^1(S)$ and $\theta^2(S)$ on the reference sphere with $S$ as referential arclength parameter. The two functions together give the spherical coordinates of a material point of the rod. Since the membrane is a fluid shell, the reference sphere used to parametrize the membrane is not necessarily the reference configuration of the membrane. This is because the in-plane fluidity precludes any canonical in-plane reference configuration. This lack of a reference configuration also causes computational issues. We adapt our earlier work, in particular the technique called \emph{gauge-fixing} \cite{dharmavaram2021gauge}, to the present problem.

We present the details of the gauge-fixed formulation and derive the weak form of the Euler-Lagrange equations, cf.\eqref{eq:weak form modified functional}. Using a combination of Loop subdivision finite element method (for the membrane) and Hermite cubic elements (for the rod), we discretize the formulation. We compute equilibrium states of using L-BFGS minimization solver. We perform a convergence test and validate the formulation on a benchmark problem \cite{das2008higher} for which semi-analytical solution is available. We also apply the method to study the behavior of a straight rod and a circular ring on the membrane. 

The focus of this work is to present a continuum framework to understand the coupled interaction between membranes and rods. In this work, we do not undertake an exhaustive study of such systems. Rather we limit ourselves to a linear isotropic constitutive model for rod. It is well known that protein filaments have a chirality and a hemitropic constitute law \cite{healey2013bifurcation, pradhan2021buckling} would be a more realistic description for rods. Adapting our formulation for such constitutive laws is straightforward nevertheless. 

Also we have limited our discussion to a system consisting of a single filament. Realistic modeling of the cortex of a cell-membrane requires a collection of cross-linked network of filaments. The cortex is a dynamic network where the junctions between the filaments are constantly changing. Modeling such systems is not difficult. One needs to add all the bending energies of the filaments and postulate some framework to add or remove crosslinkers.

We mention here some limitations of our formulation. Extending the gauge-fixing procedure to other topologies is not obvious. As discussed in \cite{dharmavaram2021gauge}, the procedure crucially relies on the existence of conformal diffeomorphism between a topologically spherical surface and the unit sphere, which is guaranteed by the Riemann-Roch theorem; see Theorem 1 \cite{dharmavaram2021gauge}. It is still an open problem if the gauge-fixing procedure can be generalized to other topologies. It is possible to circumvent the gauge-fixing procedure under certain scenarios. For example, if the deformed configuration can be described in terms of the height from some reference flat plane, i.e., the Monge representation, then a gauge-fixing procedure is not necessary. Even in spherical topologies, a radial graph representation for the deformation has been used \cite{zhao2017direct} to compute equilibrium configuration on fluid membranes, however such a description would not be able to capture asymmetrically budded shapes, e.g., Fig.~\ref{fig:bendtwist}.

In many biological systems, the cell membrane can exist in multiple fluid phases. It has long been conjectured that the in-homogeneity of lipids creates raft-like structures (\emph{raft hypothesis}) that various membrane proteins selectively bind to, thereby playing an essential role in signaling and protein transport. The boundary between the different phases can presumably bind filaments. Since the model presented above assumes that the rod freely floats on the surface of the membrane, localizing them to be at the phase boundary is not possible. If the phase boundary is known a prior, then one approach to model such a system is discussed in \cite{agrawal2009boundary}. 

\section*{Acknowledgments}
The part of the work of BLS during January 2023 was partially supported by a grant from the Simons Foundation.
BLS would like to thank the Isaac Newton Institute for Mathematical Sciences, Cambridge, for support and hospitality during January-June 2023 when a part of writing this manuscript was undertaken and was supported by EPSRC grant no EP/R014604/1.

\providecommand{\noopsort}[1]{}

\appendix

\section{Stretching and Bending Strains}
In this section, we express the stretching and bending strains explicitly in terms of the center curve of the rod (i.e. $\mathbf{r}$). The expressions we derive are used in this work to derive the first variations of these quantities and eventually used to compute the weak form in Sec.~\ref{sec:weak form}.

Recall that, since $\nu_3(S)=\parallel\mathbf{r}'(S)\parallel$, we have
\begin{equation}
\nu_3(S) = \parallel\frac{d\mathbf{r}}{dS}(S)\parallel =  \parallel{\theta^{\alpha}}'(S)\mathbf{f}_{,\alpha}(\theta^\alpha(S))\parallel,
\end{equation}
where the second equality follows from differentiating \eqref{eq:compatibility1}.

Differentiating \eqref{eq:d3=r' by nu3}, we obtain
\begin{equation}
    \mathbf{d}_3'(S) = 
    \frac{\mathbf{r}''(S)}{\nu_3}-\frac{\mathbf{r}'(S)}{\nu_3^2}\parallel\mathbf{r}'(S)\parallel',
    \label{eq:d3'}
\end{equation}
where $\parallel\mathbf{r}'(S)\parallel'=\nu_3'(S).$
Plugging in the previous equation in \eqref{eq:kappa2}, using \eqref{eq:compatibility2}, we obtain 
\begin{equation}
    \kappa_2(S) = \mathbf{d}_3'(S)\cdot \mathbf{d}_1(S)-\alpha_{31} =\frac{1}{\nu_3(S)} \mathbf{r}''(S)\cdot\mathbf{n}(\theta^\alpha(S))-\alpha_{31},
    \label{eq:kappa2 in terms of r}
\end{equation}
where we have used $\mathbf{r}'(S)\cdot \mathbf{d}_1(S)=0$ to eliminate the second term in \eqref{eq:d3'}.
Similarly, it follows by differentiating \eqref{eq:compatibility2} and using \eqref{eq:d3=r' by nu3},\eqref{eq:d2=d3xd1},\eqref{eq:kappa3} that
\begin{equation}
    \kappa_3 = \mathbf{d}_1'\cdot\mathbf{d}_2-\alpha_{12} = \theta'{}^\alpha\mathbf{n}_{,\alpha}\cdot\Big(\frac{\mathbf{r}'}{\nu_3}\times\mathbf{n}\Big)-\alpha_{12},
    \label{eq:kappa3 in terms of r}
\end{equation}
where we have used $\mathbf{d}_1'=\mathbf{n}_{,\alpha}{\theta^\alpha}
'$ (obtained by differentiating \eqref{eq:compatibility2}) in the second equality.

Finally, using (c.f.,\eqref{eq:kappa1}) $\kappa_1 = -(\mathbf{d}_3'\cdot\mathbf{d}_2-\alpha_{32})$ 
\begin{equation}
    \kappa_1 =  - \frac{1}{\nu_3}\mathbf{r}''\cdot \mathbf{d}_2 +\alpha_{32} = - \frac{1}{\nu_3}\mathbf{r}''\cdot \Big(\frac{\mathbf{r}'}{\nu_3}\times\mathbf{n}\Big) +\alpha_{32},
    \label{eq:kappa_1 in terms of r}
\end{equation}
where we have used $\mathbf{r}'(S)\cdot\mathbf{d}_2=0$ 
and \eqref{eq:compatibility2} and \eqref{eq:d2=d3xd1}.

\section{Derivation of First Variation of Strains}
\label{sec:app:deriv fv strains}
In this section, we derive explicit expressions for the variations $\delta\nu_3$, $\delta\kappa_1$, $\delta\kappa_2$, and $\delta\kappa_3$ which are used in weak form \eqref{eq:weak form modified functional}.

Note that quantities such as $\mathbf{f}$, $\bm{\eta}$, $\a_\mu$, etc., are defined on the surface $\Omega$ and are therefore functions of $(X^1,X^2)$. We recall the hat notation for evaluation of these quantities along the curve $(\theta^1(S),\theta^2(S))$ (see Sec.~\ref{ssec:redundancy}). In particular,
\begin{equation}
    \hat{\mathbf{f}}(S) = \mathbf{f}(\theta^1(S),\theta^2(S)),\quad \hat{\mathbf{f}}_{,\alpha}(S) = \mathbf{f}_{,\alpha}(\theta^1(S),\theta^2(S)),\;etc.,
\end{equation}
and, in a similar manner, for the variation $\bm{\eta}$:
\begin{equation}
    \hat{\bm{\eta}}(S) = \bm{\eta}(\theta^1(S),\theta^2(S)),\quad \hat{\bm{\eta}}_{,\alpha}(S) = \bm{\eta}_{,\alpha}(\theta^1(S),\theta^2(S)),\;etc.
\end{equation}
Note that the variation in the position of the center curve of the rod has two contributions: variation of the underlying substrate and variation in the reference configuration of the rod. That is,
\begin{equation}
    \delta\mathbf{\mathbf{r}}(S) = \hat{\bm{\eta}}(S) + \hat{\mathbf{f}}_{,\alpha}(S)\xi^\alpha(S).
\end{equation}
We drop explicit presence of $S$ in the following while writing various expressions. Differentiating the previous equation, we obtain
\begin{equation}
    (\delta\mathbf{r})'= \Big(\bm{\hat\eta}_{,\alpha}\theta^\alpha{}'+\hat{\mathbf{f}}_{,\alpha\beta}\theta^\beta{}'\xi^\alpha+\hat{\mathbf{f}}_{,\alpha}\xi^\alpha{}'\Big),
\end{equation}
and differentiating it again, we get
\begin{equation}
   (\delta\mathbf{r})'' =\Big(\hat{\bm{\eta}}_{,\alpha\beta}\theta'{}^\alpha\theta^\beta{}' + \hat{\bm{\eta}}_{,\alpha}\theta''{}^\alpha + \hat{\mathbf{f}}_{,\alpha\beta\gamma}\theta'{}^\beta\theta'{}^\gamma\xi^\alpha+\hat{\mathbf{f}}_{,\alpha\beta}\theta''{}^\beta\xi^\alpha + 2\hat{\mathbf{f}}_{,\alpha\beta}\theta'{}^\beta\xi^\alpha{}'+ \hat{\mathbf{f}}_{,\alpha}\xi^\alpha{}''\Big).
\label{eq:delta r''}
\end{equation}
The variations in the directors are as follows:
\begin{align}
    \delta\mathbf{d}_1&=\delta\hat{\mathbf{n}} = \delta_s\hat{\mathbf{n}} + \hat{\mathbf{n}}_{,\alpha}\xi^\alpha = -(\hat{\a}^\mu\otimes\hat{\mathbf{n}})\cdot\hat{\bm{\eta}}_{,\mu} + \hat{\mathbf{n}}_{,\alpha}\xi^\alpha,
    \label{eq:del d1}\\
    \delta\mathbf{d}_2 &= \delta(\mathbf{d}_3\times\hat{\mathbf{n}})=\delta\mathbf{d}_3\times\hat{\mathbf{n}} + \mathbf{d}_3\times\delta\hat{\mathbf{n}},
     \label{eq:del d2}\\
    \delta\mathbf{d}_3 &= \frac{\delta\mathbf{r}'}{\nu_3}-\frac{\mathbf{r}'}{\nu_3^2}\delta\nu_3,
     \label{eq:del d3}
\end{align}
where we have used the surface variation of the normal, $\delta_s\hat{\mathbf{n}}:=-(\hat{\a}^\mu\otimes\hat{\mathbf{n}})\cdot\hat{\bm{\eta}}_{,\mu}$, in \eqref{eq:del d1} (see (27) of \cite{dharmavaram2021gauge} for details). 

We note the following results for the variations.
\begin{prop}
\begin{equation}
\delta\nu_3=\frac{{\mathbf{r}}'}{\nu_3}\cdot \Big(\hat{\bm{\eta}}_{,\alpha}{\theta^\alpha}'+\hat{\mathbf{f}}_{,\alpha\beta}{\theta^\beta}'\xi^\alpha+\hat{\mathbf{f}}_{,\alpha}{\xi^\alpha}'\Big)
\label{eq:delta nu3}
\end{equation}
\label{prop:nu3}
\end{prop}
\begin{proof}
This directly follows from taking the variation of $\nu_3=\sqrt{\mathbf{r}'\cdot\mathbf{r}'}$. That is,
$$
\delta\nu_3 = \delta\sqrt{\mathbf{r}'\cdot\mathbf{r}'} =\frac{1}{\parallel{\mathbf{r}}'\parallel} {\mathbf{r}}'\cdot\delta{\mathbf{r}}' 
= \frac{{\mathbf{r}}'}{\parallel{\mathbf{r}'}\parallel}\cdot \Big(\bm{\eta}(\theta^\alpha(\varphi))+\hat{\mathbf{f}}_{,\alpha}\xi^\alpha\Big)',
$$
where \eqref{eq:delta nu3} follows from using $\nu_3=\parallel\mathbf{r}'\parallel$ and by differentiating the terms in the round brackets.
\end{proof}

\begin{prop}
\begin{equation}
    \delta\kappa_2 = \frac{1}{\nu_3}\Big\{(\delta\mathbf{r})''\cdot\hat{\mathbf{n}}+\mathbf{r}''\cdot\delta\hat{\mathbf{n}}-(\kappa_2+\alpha_{31})\delta\nu_3\Big\},
\end{equation}
where $\delta\nu_3$ is given by \eqref{eq:delta nu3}, $(\delta\mathbf{r})''$ and $\delta\hat{\mathbf{n}}$ are given by \eqref{eq:delta r''} and \eqref{eq:del d1}, respectively.
\end{prop}
\label{prop:k2}
\begin{proof}
Taking the variation of \eqref{eq:kappa2 in terms of r}, we obtain
\begin{equation}
    \delta\kappa_2 = \frac{1}{\nu_3}\Big((\delta\mathbf{r})''\cdot\hat{\mathbf{n}} + \mathbf{r}''\cdot\delta\hat{\mathbf{n}}\Big)-\mathbf{r}''\cdot\hat{\mathbf{n}}\frac{1}{\nu_3^2}\delta\nu_3.
    \label{eq:dkappa2temp}
\end{equation}
Using the expression for $\kappa_{2}$ as given in \eqref{eq:kappa2 in terms of r} into the last term of \eqref{eq:dkappa2temp}, we obtain
\begin{equation}
    \delta\kappa_2 = \frac{1}{\nu_3}\Big\{(\delta\mathbf{r})''\cdot\hat{\mathbf{n}}+\mathbf{r}''\cdot\delta\hat{\mathbf{n}}-(\kappa_2+\alpha_{31})\delta\nu_3\Big\}.
\end{equation}
\end{proof}
\begin{prop}
\begin{equation}
    \delta\kappa_1 = -\frac{1}{\nu_3}\Big\{(\delta\mathbf{r})''\cdot\mathbf{d}_2+\mathbf{r}''\cdot\delta\mathbf{d}_2+(\kappa_1-\alpha_{32})\delta\nu_3\Big\},
\end{equation}
\end{prop}
\label{prop:k1}
\begin{proof}
Taking the variation of \eqref{eq:kappa_1 in terms of r}, we obtain
\begin{equation}
    \delta\kappa_1 = -\frac{1}{\nu_3}\Big((\delta\mathbf{r})''\cdot\mathbf{d}_2 + \mathbf{r}''\cdot\delta\mathbf{d}_2\Big)+\mathbf{r}''\cdot\mathbf{d}_2\frac{1}{\nu_3^2}\delta\nu_3.
    \label{eq:dkappa1temp}
\end{equation}
Using the expression for $\kappa_1$ as given in \eqref{eq:kappa_1 in terms of r} into the last term of \eqref{eq:dkappa1temp}, we obtain
\begin{equation}
    \delta\kappa_1 = -\frac{1}{\nu_3}\Big\{(\delta\mathbf{r})''\cdot\mathbf{d}_2+\mathbf{r}''\cdot\delta\mathbf{d}_2+(\kappa_1-\alpha_{32})\delta\nu_3\Big\}.
\end{equation}

\end{proof}
\begin{prop}
\begin{equation}
       \delta\hat{\mathbf{n}}_{,\alpha} = -(\hat{\mathbf{f}}^\mu_{,\alpha}\otimes \hat{\mathbf{n}} + \hat{\mathbf{f}}^\mu\otimes\hat{\mathbf{n}}_{,\alpha})\cdot \hat{\bm{\eta}}_{,\mu} - (\hat{\mathbf{f}}^\mu\otimes \hat{\mathbf{n}})\cdot\hat{\bm{\eta}}_{,\mu\alpha} + \hat{\mathbf{n}}_{,\alpha\beta}\xi^\beta,
\end{equation}
\label{prop4}
\end{prop}
\begin{proof}
Recall from \eqref{eq:delta_d_alpha} that 
\begin{equation}
  \delta\mathbf{n}_{,\alpha} = -(\mathbf{a}^\beta_{,\alpha}\otimes \mathbf{n} + \mathbf{a}^\beta\otimes\mathbf{n}_{,\alpha})\cdot\bm{\eta}_{,\beta} - (\mathbf{a}^\beta\otimes\mathbf{n})\cdot\bm{\eta}_{,\beta\alpha}.
  \label{eq:delta_d_alpha2}
\end{equation}
It then follows that 
\begin{equation}
     \delta\hat{\mathbf{n}}_{,\alpha} = -(\hat{\mathbf{f}}^\mu_{,\alpha}\otimes \hat{\mathbf{n}} + \hat{\mathbf{f}}^\mu\otimes\hat{\mathbf{n}}_{,\alpha})\cdot \hat{\bm{\eta}}_{,\mu} - (\hat{\mathbf{f}}^\mu\otimes \hat{\mathbf{n}})\cdot\hat{\bm{\eta}}_{,\mu\alpha} + \hat{\mathbf{n}}_{,\alpha\beta} \xi^\beta,
\end{equation}
where the last term arises due to the variation in $\bm{\theta}$.
\end{proof}
\begin{prop}
\begin{equation}
    \delta\kappa_3 = \mathbf{d}_2\cdot\delta\hat{\mathbf{n}}_{,\alpha}\theta^{\alpha}{}' + \mathbf{d}_2\cdot\hat{\mathbf{n}}_{,\alpha}\xi^{\alpha}{}' + \theta^{\alpha}{}'\hat{\mathbf{n}}_{,\alpha}\cdot\delta\mathbf{d}_2,
\end{equation}
where $\delta\hat{\mathbf{n}}_{,\alpha}$ is computed according to Prop.~\ref{prop4}.
\label{prop:k3}
\end{prop}
\begin{proof}
Taking the variation of \eqref{eq:kappa3 in terms of r}, we obtain
\begin{equation}
\delta\kappa_3 = \delta(\mathbf{d}_1')\cdot\mathbf{d}_2 + \mathbf{d}_1'\cdot\delta\mathbf{d}_2 = \delta(\hat{\mathbf{n}}_{,\alpha}{\theta^\alpha}')\cdot\mathbf{d}_2 + \delta\mathbf{d}_2\cdot\hat{\mathbf{n}}_{,\alpha}\theta^{\alpha}{}',
\end{equation}
where we have used the identity $\mathbf{d}_1'=\hat{\mathbf{n}}_{,\alpha}{\theta^\alpha}
'$ to obtain the second inequality. Next, expanding the variation of the term in the round bracket, we obtain
\begin{equation}
    \delta\kappa_3 = (\delta\hat{\mathbf{n}}_{,\alpha}\theta^{\alpha}{}'(S) + \hat{\mathbf{n}}_{,\alpha}\xi^{\alpha}{}')\cdot\mathbf{d}_2 + \hat{\mathbf{n}}_{,\alpha}\theta^{\alpha}{}'(S)\cdot\delta\mathbf{d}_2
\end{equation}
\end{proof}
\end{document}